# Secondary Analytic Indices

John Lott[*]
Department of Mathematics
University of Michigan
Ann Arbor, MI 48109-1003
USA
lott@math.lsa.umich.edu

**Abstract**

We define real-valued characteristic classes of flat complex vector bundles and flat real vector bundles with a duality structure. We construct pushforwards of such vector bundles with vanishing characteristic classes. These pushforwards involve the analytic torsion form in the first case and the eta-form of the signature operator in the second case. We show that the pushforwards are independent of the geometric choices made in the constructions and hence are topological in nature. We give evidence that in the first case, the pushforwards are given topologically by the Becker-Gottlieb-Dold transfer.

## 1 Introduction

In this paper we construct analytic indices for secondary, as opposed to primary, index theorems. An example of a primary index theorem is the Atiyah-Singer families index theorem [3]. Roughly speaking, when both sides of a primary index theorem vanish over the reals, a secondary index theorem gives more refined information. The construction of the secondary analytic index involves so-called secondary invariants, such as the Bismut-Cheeger

[*]Partially supported by NSF grant DMS-9403652.

eta-form [7] or the analytic torsion form of the author and Bismut [9]. We construct secondary analytic indices in two interesting cases. We make a literal analogy between the columns of the following table :

Table 1

Primary Data

| **C** − vector bundle | flat **C** − vector bundle | flat **R** − duality bundle |
|---|---|---|
| topological K theory | algebraic K theory | ? |
| Chern character | Borel classes | $p$ character |
| Dirac type operator $D$ | de Rham operator | signature operator $\sigma$ |
| spin$^c$ | vacuous | oriented |
| $\widehat{A}$ class | Euler class | $L$ class |
| AS families index thm. | index thm. of [9] | families index thm. for $\sigma$ |
| Gysin map | transfer of [4, 15] | Gysin map |

Secondary Data

| $K^{-1}(\,\cdot\,;\mathbf{R}/\mathbf{Z})$ | $\overline{K}^0_R$ | $\overline{L}^0_\epsilon$ |
|---|---|---|
| reduced eta invariant | analytic torsion | eta invariant of $\sigma$ |
| eta form of $D$ | analytic torsion form | eta form of $\sigma$ |

To explain the meaning of this table, consider the first column. The primary data is explained, for example, in the book of Berline-Getzler-Vergne[5]. Suppose that $M \xrightarrow{\pi} B$ is a smooth fiber bundle with even-dimensional connected closed fibers $Z$. Suppose in addition that the vertical tangent bundle $TZ$ has a spin$^c$-structure. Then there is an analytically-defined pushforward $\pi_! : K^0(M) \to K^0(B)$. The Atiyah-Singer families index theorem says that $\pi_!$ is the same as the topologically-defined Gysin map $G : K^0(M) \to K^0(B)$. As a consequence, if $E$ is a virtual vector bundle on $M$ then one obtains an identity in H$^*(B;\mathbf{R})$ :

$$\mathrm{ch}(\pi_!(E)) = \int_Z \widehat{A}(TZ) \wedge \mathrm{ch}(E). \tag{1}$$

If ch($E$) vanishes in H$^*(M; \mathbf{R})$ then both sides of (1) vanish. This suggests focusing on (virtual) vector bundles with vanishing Chern character. Karoubi showed that they are intimately related to K-theory with coefficients. In fact, there is a geometric description of $K^{-1}(M; \mathbf{R}/\mathbf{Z})$ in terms of such vector bundles. Its generators are given by (virtual) Hermitian vector bundles on $M$ with connection whose Chern characters are written explicitly as exact forms [20, 22].

In [22] we showed that one can detect elements of $K^{-1}(M; \mathbf{R}/\mathbf{Z})$ analytically by mapping a closed odd-dimensional spin$^c$-manifold $X$ into $M$ and computing Atiyah-Patodi-Singer eta-invariants of Dirac-type operators on $X$, reduced mod $\mathbf{Z}$. Furthermore, if $M \xrightarrow{\pi} B$ is a fiber bundle as above then we constructed an analytic pushforward $\pi_! : K^{-1}(M; \mathbf{R}/\mathbf{Z}) \to K^{-1}(B; \mathbf{R}/\mathbf{Z})$. The definition of this pushforward involved the eta-form of the fiberwise Dirac-type operators. The secondary index theorem of [22] says that $\pi_!$ equals the Gysin map $G : K^{-1}(M; \mathbf{R}/\mathbf{Z}) \to K^{-1}(B; \mathbf{R}/\mathbf{Z})$.

Thus the first column of Table 1 illustrates a method which starts with a primary index theorem and produces a secondary index theorem. In this paper we apply this method to two other primary index theorems. The first is an index theorem for flat complex vector bundles which appeared in [9] and is summarized in Appendix A of the present paper. This gives rise to the second column of Table 1. The relevant characteristic classes of flat complex vector bundles are known in algebraic K-theory as the Borel classes. Given a smooth fiber bundle $M \xrightarrow{\pi} B$ with connected closed fibers $Z$ and a flat complex vector bundle $E$ on $M$, the "primary" analytic pushforward of $E$ is simply given by the cohomology groups of the fibers, with value in the fiberwise restrictions of $E$. These cohomology groups form flat complex vector bundles on $B$.

To produce a secondary index, in Section 2 we consider a K-theory generated by flat complex vector bundles $E$ on a manifold $M$ whose Borel classes are written explicitly as exact forms. Actually, in order to produce interesting examples we must assume that the flat vector bundles have a more rigid structure. As part of the data, we assume that there is a ring $R$, a representation $\rho : R \to \text{End}(\mathbf{C}^n)$ and a local system $F$ of $R$-modules on $M$ such that $E = F \otimes_\rho \mathbf{C}^n$. We then define a K-group $\overline{K}_R^0(M)$ of such objects with trivialized Borel classes. Our work is similar in spirit to that of Gillet and Soulé in the holomorphic case [33], with their arithmeticity assumption being replaced by our use of the ring $R$.

We then define the secondary analytic pushforward. Given a fiber bundle $M \xrightarrow{\pi} B$ as just above, we add some additional geometric structure in the form of a horizontal distribution $T^H M$ and a family of vertical Riemannian metrics $g^{TZ}$. We then use the analytic torsion form of [9] to define an analytic pushforward $\pi_! : \overline{K}^0_R(M) \to \overline{K}^0_R(B)$. We show that $\pi_!$ is independent of the choices of $T^H M$ and $g^{TZ}$ and hence is topological in nature. We conjecture that $\pi_!$ is related to the Becker-Gottlieb-Dold transfer [4, 15] in a certain generalized cohomology theory which is related to algebraic K-theory. We show that this conjecture is true when $B$ is a point. This uses the fact that one knows what the analytic torsion form is in this case, namely the Reidemeister torsion.

Section 3 deals with the third column of Table 1. We consider a primary index theory which is based on flat real vector bundles with a Poincaré-duality-type structure. We first define an analog of the Chern character for such bundles. Let $M \xrightarrow{\pi} B$ be a smooth fiber bundle with even-dimensional connected closed fibers $Z$ such that the vertical tangent bundle $TZ$ is oriented. Given a flat "duality" vector bundle $E$ on $M$, the analytic pushforward of $E$ is again constructed from the cohomology groups of the fibers, with value in the fiberwise restrictions of $E$; this uses the fact that the orientation of the fibers gives a Poincaré-duality structure on the cohomology groups. The primary index theorem is simply the Atiyah-Singer families index theorem when applied to the vertical signature operators. It was previously considered by Atiyah in the special case when $E$ is trivial, in order to explain the nonmultiplicativity of the signature of fiber bundles [1]. (The pushforward of "duality" vector bundles which are fiberwise flat but not globally flat was considered by Lusztig for applications to the Novikov conjecture [26], but this differs from what we do.)

In order to define a secondary index, we consider a group $\overline{L}^0_\epsilon(M)$ generated by flat duality vector bundles on $M$ whose Chern-type characters are written explicitly as exact forms. We show that one can detect elements of $\overline{L}^0_\epsilon(M)$ by mapping a closed odd-dimensional oriented manifold $X$ into $M$ and computing eta-invariants of twisted tangential signature operators on $X$. This perhaps gives some hint to the topological meaning of real-valued eta-invariants. Given a fiber bundle $M \xrightarrow{\pi} B$ as just above, we add a horizontal distribution $T^M$ and a family of vertical Riemannian metrics $g^{TZ}$. We then use the eta-form of the vertical signature operators to define an analytic pushforward $\pi_! : \overline{L}^0_\epsilon(M) \to \overline{L}^0_{\epsilon\epsilon_n}(B)$. We show that $\pi_!$ is independent of $T^H M$ and $g^{TZ}$. However, it is not clear to us what the purely topological equivalent of $\pi_!$ should be.

The detailed descriptions of the subsections of this paper appear at the beginnings of Sections 2 and 3. The contents are:

1. Introduction

2. Flat Complex Vector Bundles
2.1 $\overline{K}$-Groups
2.2 Analytic Pushforward of $\overline{K}$-Groups
2.3 Possible Topological Interpretation

3. Flat Duality Bundles
3.1 Characteristic Classes of Flat Duality Bundles
3.2 Flat Duality Superconnections
3.3 $\overline{L}$-Groups and Real-Valued Eta-Invariants
3.4 Number Operators
3.5 Fiber Bundles
3.6 Analytic Pushforward of $\overline{L}$-Groups

A. Results from [9]
A.1 Characteristic Classes of Flat Complex Vector Bundles
A.2 The Superconnection Formalism
A.3 Characteristic Classes and Torsion Forms of Flat Superconnections
A.4 Fiber Bundles

I thank Wolfgang Lück, Christophe Soulé, Michael Weiss and Bruce Williams for helpful discussions. It should be clear that many of the constructions in this paper are inspired by ideas of Gillet and Soulé. I thank Jean-Pierre Bourguignon and the IHES for their hospitality while part of this research was performed.

## 2 Flat Complex Vector Bundles

Let $R$ be a ring, satisfying certain conditions to be specified. Let $\rho : R \to \text{End}(\mathbf{C}^n)$ be a representation of $R$. We consider local systems $F$ of $R$-modules on a manifold $M$ and their complexifications $F_{\mathbf{C}} = F \otimes_\rho \mathbf{C}^n$. Given a Hermitian metric $h^F$ on the flat complex vector bundle $F_{\mathbf{C}}$, certain characteristic forms $c_*(\nabla^F, h^F) \in \Omega^{odd}(M)$ were defined in [9] (see Appendix A.1). In Subsection 2.1 we define a group $\overline{K}^0_R(M)$ which essentially consists of pairs $(F, h^F)$ along with an explicit writing of $c_*(\nabla^F, h^F)$ as an exact form. In fact, we first define a group $\widehat{K}^0_R(M)$ and a map $c' : \widehat{K}^0_R(M) \to \Omega^{odd}(M)$, and then put $\overline{K}^0_R(M) = \text{Ker}(c')$. We show that $\overline{K}^0_R(M)$ is a homotopy-invariant of $M$ and compute $\overline{K}^0_R(\text{pt.})$. In Subsection 2.2 we start with a smooth fiber bundle $M \xrightarrow{\pi} B$ with connected closed fibers $Z$. Given the additional geometric data of a horizontal distribution $T^H M$ on $M$ and vertical Riemannian metrics $g^{TZ}$ on the fibers, we construct an analytic pushforward $\pi_! : \overline{K}^0_R(M) \to \overline{K}^0_R(B)$. This analytic pushforward involves the analytic torsion form of the fiber bundle. We show that $\pi_!$ is independent of $T^H M$ and $g^{TZ}$ and hence depends only on the smooth topological structure of the fiber bundle. Using the Cheeger-Müller theorem, we compute $\pi_!$ when $B$ is a point. In Subsection 2.3 we give a plausible topological equivalent of $\pi_!$ in terms of the Becker-Gottlieb-Dold transfer.

### 2.1 $\overline{K}$-Groups

Let $R$ be a right-Noetherian ring. Suppose that $R$ is right-regular, meaning that every finitely-generated right-$R$-module has a finite resolution by finitely-generated projective right-$R$-modules. (A relevant example is $R = \mathbf{Z}$.) Hereafter, all right-$R$-modules will be taken to be finitely-generated.

Let $M$ be a connected smooth manifold and let $F$ be a local system on $M$ modeled on a right-$R$-module $V$ [34, p. 58]. Let $m_0$ be a basepoint in $M$, put $\Gamma = \pi_1(M, m_0)$ and let $\widetilde{M}$ be the universal cover of $M$. If $V$ is the fiber of $F$ over $m_0$ then $V$ is a right-$R$-module and a left-$\mathbf{Z}\Gamma$-module, and $F = \widetilde{M} \times_\Gamma V$. Conversely, given an $R\Gamma$-module $V$ which is finitely-generated over $R$, we can construct a local system $F$ of right-$R$-modules on $M$ by $F = \widetilde{M} \times_\Gamma V$.

**Definition 1** *The group $K^0_R(M)$ is the quotient of the free abelian group generated by local systems of right-$R$-modules on $M$, by the subgroup gener-* ated by the relations $F^2 - F^1 - F^3$ whenever there is a short exact sequence of local systems

$$0 \longrightarrow F^1 \longrightarrow F^2 \longrightarrow F^3 \longrightarrow 0. \tag{2}$$

We can identify $K_R^0(M)$ as the Grothendieck group of isomorphism classes of $R\Gamma$-modules which are finitely-generated over $R$. With our regularity assumption on $R$, we have $K_R^0(\text{pt.}) \cong K_0(R)$.

If $R$ and $R'$ are two rings obeying our assumptions, there is a natural product

$$K_R^0(M) \times K_{R'}^0(M) \longrightarrow K_{R \otimes R'}^0(M) \tag{3}$$

generated by

$$[F] \times [F'] = [F \otimes F']. \tag{4}$$

Let $\rho : R \to \text{End}(\mathbf{C}^n)$ be a ring homomorphism such that $\mathbf{C}^n$ is flat when considered as an $R$-module. (A relevant example is when $R = \mathbf{Z}$ and $\rho : \mathbf{Z} \to \text{End}(\mathbf{C})$ is the standard inclusion.) If $V$ is a right-$R$-module, let $V_\mathbf{C}$ denote the complex vector space $V_\mathbf{C} = V \otimes_\rho \mathbf{C}^n$. Let $F_\mathbf{C} = F \otimes_\rho \mathbf{C}^n$ denote the complex vector bundle on $M$ associated to the local system $F$. It has a flat connection $\nabla^F$. Define the characteristic class $c(\nabla^F) \in \text{H}^{odd}(M; \mathbf{R})$ as in Definition 38 of the Appendix. Given a short exact sequence (2), by tensoring over $\rho$ we obtain a short exact sequence of flat complex vector bundles on $M$

$$0 \longrightarrow F_\mathbf{C}^1 \xrightarrow{v} F_\mathbf{C}^2 \xrightarrow{v} F_\mathbf{C}^3 \longrightarrow 0. \tag{5}$$

**Proposition 1** *The assignment of $c(\nabla^F)$ to $F$ extends to a map $c : K_R^0(M) \to \text{H}^{odd}(M; \mathbf{R})$.*

**Pf.** : Given a sequence (2), we must show that $c(\nabla^{F^2}) = c(\nabla^{F^1}) + c(\nabla^{F^3})$. This follows from (294). ∎

Let $h^F$ be a positive-definite Hermitian metric on $F_\mathbf{C}$. Define $c(\nabla^F, h^F) \in \Omega^{odd}(M)$ as in Definition 37.

**Definition 2** *A $\widehat{K}_R^0(M)$-generator is a triple $\mathcal{F} = (F, h^F, \eta)$ where*

- *$F$ is a local system of right-$R$-modules on $M$.*

- *$h^F$ is a positive-definite Hermitian metric on $F_\mathbf{C}$.*

- $\eta$ is an element of $\Omega^{even}(M)/\operatorname{im}(d)$.

Given a short exact sequence (5), choose Hermitian metrics $\left\{h^{F^i}\right\}_{i=1}^3$ on $\{F_{\mathbf{C}}^i\}_{i=1}^3$. Define the torsion form $T_f(A', h^F) \in \Omega^{even}(M)$ as in Definition 46. Its salient property is that it satisfies equation (294).

**Definition 3** *A $\widehat{K}_R^0(M)$-relation is given by three $\widehat{K}_R^0(M)$-generators $\mathcal{F}^1$, $\mathcal{F}^2$ and $\mathcal{F}^3$, along with a short exact sequence (2) such that*

$$\eta_2 = \eta_1 + \eta_3 + T_f(A', h^F). \tag{6}$$

**Definition 4** *The group $\widehat{K}_R^0(M)$ is the quotient of the free abelian group generated by the $\widehat{K}_R^0(M)$-generators, by the subgroup generated by the $\widehat{K}_R^0(M)$-relations $\mathcal{F}^2 - \mathcal{F}^1 - \mathcal{F}^3$.*

The reader can compare the definition of $\widehat{K}_R^0(M)$ with that given in [17, Defn. 6.1] for the analogous holomorphic case.

**Proposition 2** *The assignment of $c(\nabla^F, h^F) - d\eta$ to $(F, h^F, \eta)$ extends to a map $c' : \widehat{K}_R^0(M) \to \Omega^{odd}(M)$.*

**Pf. :** Given a sequence (2) satisfying (6), we must show that

$$c(\nabla^{F^2}, h^{F^2}) - c(\nabla^{F^1}, h^{F^1}) - c(\nabla^{F^3}, h^{F^3}) = d\eta_2 - d\eta_1 - d\eta_3. \tag{7}$$

This follows from (294). ∎

**Definition 5** *Let $\overline{K}_R^0(M)$ be the kernel of $c'$.*

There is a complex

$$\mathrm{H}^{even}(M; \mathbf{R}) \xrightarrow{a} \overline{K}_R^0(M) \xrightarrow{b} K_R^0(M) \xrightarrow{c} \mathrm{H}^{odd}(M; \mathbf{R}), \tag{8}$$

where $a(\sigma) = [(0, 0, \sigma)]$ and $b\left(\sum_i n_i \left[(F^i, h^{F^i}, \eta^i)\right]\right) = \sum_i n_i [F^i]$.

**Proposition 3** *The complex (8) is exact.*

**Pf. :** Exactnesss at $K_R^0(M)$ : Given $\sum_i n_i[F^i] \in K_R^0(M)$ such that $c(\sum_i n_i[F^i]) = 0$, choose arbitrary Hermitian metrics $\{h^{F^i}\}$ on the $\{(F^i)_\mathbf{C}\}$. As $\sum_i n_i c(\nabla^{F^i}, h^{F^i}) \in \Omega^{odd}(M)$ represents $c(\sum_i n_i[F^i])$ in de Rham cohomology, there is a $\tau \in \Omega^{even}(M)/\mathrm{im}(d)$ such that $\sum_i n_i c(\nabla^{F^i}, h^{F^i}) = d\tau$. Then $[(0, 0, \tau)] + \sum_i n_i[(F^i, h^{F^i}, 0)]$ is an element of $\overline{K}_R^0(M)$ whose image under $b$ is $\sum_i n_i[F^i]$.

Exactness at $\overline{K}_R^0(M)$ : Given $\sum_i n_i[(F^i, h^{F^i}, \eta^i)] \in \overline{K}_R^0(M)$ such that $\sum_i n_i[F^i]$ vanishes in $K_R^0(M)$, there are local systems $\{G^j, H^j, I^j\}$ on $M$ and short exact sequences
$$0 \longrightarrow G^j \longrightarrow H^j \longrightarrow I^j \longrightarrow 0 \tag{9}$$
so that
$$\sum_i n_i F^i = \sum_j m_j \left(G^j - H^j + I^j\right) \tag{10}$$
for some integers $\{m_j\}$. We also have
$$\sum_i n_i c(\nabla^{F^i}, h^{F^i}) = d \sum_i n_i \eta^i \tag{11}$$

Put Hermitian metrics $\{h^{G^j}, h^{H^j}, h^{I^j}\}$ on $\{(G^j)_\mathbf{C}, (H^j)_\mathbf{C}, (I^j)_\mathbf{C}\}$ so that (with a slight abuse of terminology) coincident terms in (10) have the same Hermitian metric. Then
$$\sum_i n_i c(\nabla^{F^i}, h^{F^i}) = \sum_j m_j \left(c(\nabla^{G^j}, h^{G^j}) - c(\nabla^{H^j}, h^{H^j}) - c(\nabla^{I^j}, h^{I^j})\right). \tag{12}$$

Let $T^j \in \Omega^{even}(B)/\mathrm{im}(d)$ be the torsion form of the sequence (9). By (294),
$$dT^j = -c(\nabla^{G^j}, h^{G^j}) + c(\nabla^{H^j}, h^{H^j}) - c(\nabla^{I^j}, h^{I^j}). \tag{13}$$

From (11), (12) and (13), we have that $\tau \equiv \sum_i n_i \eta^i + \sum_j m_j T^j$ lies in $\mathrm{H}^{even}(B; \mathbf{R})$. Then in $\overline{K}_R^0(M)$, we have
$$\begin{aligned}\sum_i n_i[(F^i, h^{F^i}, \eta^i)] &= a(\tau) + \sum_j m_j \left([(G^j, h^{G^j}, 0)] - [(H^j, h^{H^j}, T^j)] \right. \\ &\qquad\qquad\qquad\qquad\qquad \left. + [(I^j, h^{I^j}, 0)]\right) \\ &= a(\tau). \end{aligned} \tag{14}$$

The proposition follows. ∎

We now show the homotopy invariance of $\overline{K}_R^0(\cdot)$.

**Proposition 4** *Define $i_0 : M \to [0,1] \times M$ by $i_0(m) = (0, m)$ and $i_1 : M \to [0,1] \times M$ by $i_1(m) = (1, m)$. Then for any $\widetilde{k} \in \overline{K}_R^0([0,1] \times M)$, one has $i_0^* \widetilde{k} = i_1^* \widetilde{k}$.*

**Pf.** : Letting $s$ be a coordinate on $[0,1]$, we can write exterior differentiation on $[0,1] \times M$ as
$$\widetilde{d} = ds \wedge \partial_s + d. \tag{15}$$
Let us write $\widetilde{k} = \sum_j n_j [(\widetilde{F^j}, h^{\widetilde{F^j}}, \widetilde{\eta}^j)]$, with
$$\sum_j n_j \, c(\nabla^{\widetilde{F^j}}, h^{\widetilde{F^j}}) = \widetilde{d} \sum_j n_j \widetilde{\eta}^j. \tag{16}$$

Denote the restriction of $[(\widetilde{F^j}, h^{\widetilde{F^j}}, \widetilde{\eta}^j)]$ to $\{s\} \times M$ by $[(F^j, h_s^{F^j}, \eta_s^j)]$. Then equations (15) and (16) imply that modulo $\operatorname{im}(d)$,
$$\frac{\partial}{\partial s} \sum_j n_j \eta_s^j = i_{\partial_s} \sum_j n_j \, c(\nabla^{\widetilde{F^j}}, h^{\widetilde{F^j}}). \tag{17}$$

Let $T_f(A', h_s^{F^j})$ be the torsion form of the sequence
$$0 \to (F^j, h_0^{F^j}) \xrightarrow{Id} (F^j, h_s^{F^j}) \to 0 \to 0. \tag{18}$$

One can show from (294) that
$$\frac{\partial}{\partial s} \sum_j n_j \, T_f(A', h_s^{F^j}) = i_{\partial_s} \sum_j n_j \, c(\nabla^{\widetilde{F^j}}, h^{\widetilde{F^j}}). \tag{19}$$

Thus for all $s \in [0, 1]$,
$$\sum_j n_j \, \eta_s^j = \sum_j n_j \, \eta_0^j + \sum_j n_j \, T_f(A', h_s^{F^j}). \tag{20}$$

It follows that in $\widehat{K}_R^0(M)$,
$$\sum_j n_j [(F^j, h_s^{F^j}, \eta_s^j)] = \sum_j n_j [(F^j, h_0^{F^j}, \eta_0^j)]. \tag{21}$$

The proposition follows. ∎

**Corollary 1** *Let $Z$ be a smooth connected manifold. Given a smooth map $\phi : Z \to M$ and an element $k \in \overline{K}_R^0(M)$, the pullback $\phi^* k \in \overline{K}_R^0(Z)$ only depends on the (smooth) homotopy class of $\phi$.*

**Corollary 2** *If $f : M \to M'$ is a (smooth) homotopy equivalence then it induces an isomorphism $f^* : \overline{K}_R^0(M') \to \overline{K}_R^0(M)$.*

Let $R$ and $R'$ be two rings obeying our assumptions. Then there is a natural product
$$\widehat{K}_R^0(M) \times \widehat{K}_{R'}^0(M) \longrightarrow \widehat{K}_{R \otimes R'}^0(M) \tag{22}$$
generated by
$$[(F, h^F, \eta)] \times [(F', h^{F'}, \eta')] = [(F \otimes F', h^F \otimes h^{F'}, \mathrm{rk}(F_{\mathbf{C}}) \, \eta' + \mathrm{rk}(F'_{\mathbf{C}}) \, \eta)]. \tag{23}$$
This passes to a product
$$\overline{K}_R^0(M) \times \overline{K}_{R'}^0(M) \longrightarrow \overline{K}_{R \otimes R'}^0(M). \tag{24}$$

We now consider the special case when $M$ is a point.

**Definition 6** *Let $K_0^{vol}(R)$ be the Grothendieck group of isomorphism classes of pairs $(V, \mathrm{vol})$ consisting of a right-$R$-module $V$ and a volume form $\mathrm{vol}$ on $V_{\mathbf{C}}$.*

**Proposition 5** $\overline{K}_R^0(\mathrm{pt.})$ *is isomorphic to* $K_0^{vol}(R)$.

**Pf.** : A generator for $\overline{K}_R^0(\mathrm{pt.})$ is a triple $(V, h^V, \eta)$ where $V$ is a right-$R$-module, $h^V$ is a Hermitian metric on $V_{\mathbf{C}}$ and $\eta \in \mathbf{R}$. With our regularity assumption on $R$, we may assume that $V$ is projective. Let $\mathrm{vol}(h^V)$ be the induced volume form on $V_{\mathbf{C}}$. To $(V, h^V, \eta)$ we assign the pair $(V, e^{-\eta} \mathrm{vol}(h^V))$. We claim that this passes to a map from $\overline{K}_R^0(\mathrm{pt.})$ to $K_0^{vol}(R)$. To see this, suppose that we have a short exact sequence (2) satisfying (6). In this case, (6) becomes
$$\eta_2 = \eta_1 + \eta_3 + \ln\left(\mathrm{vol}(h^{F^2})/\mathrm{vol}(h^{F^1})\mathrm{vol}(h^{F^3})\right), \tag{25}$$
where we use an obvious notation. Thus we get a relationship of pairs
$$0 \to (F^1, e^{-\eta_1}\mathrm{vol}(h^{F^1})) \to (F^2, e^{-\eta_2}\mathrm{vol}(h^{F^2})) \to (F^3, e^{-\eta_3}\mathrm{vol}(h^{F^3})) \to 0, \tag{26}$$

which is what was needed to be shown.

It is now straightforward to see that this gives the isomorphism of the proposition. ∎

The ring homomorphism $\rho$ induces a map $\rho_* : K_1(R) \to K_1(\mathbf{C}) \cong \mathbf{C}^*$.

**Proposition 6** *There is an exact sequence*

$$K_1(R) \xrightarrow{\ln |\rho_*|} \mathbf{R} \xrightarrow{a} \overline{K}^0_R(\text{pt.}) \xrightarrow{b} K_0(R) \xrightarrow{c} 0. \tag{27}$$

**Pf.** : It remains to show exactness at $\mathbf{R}$. Recall that $K_1(R)$ is generated by automorphisms $A$ of projective right-$R$-modules $V$. If vol is a volume form on $V_{\mathbf{C}}$ then $A^*\text{vol} = |\rho_*(A)|$ vol. Thus $A$ gives an isomorphism between the pairs $(V, \text{vol})$ and $(V, |\rho_*(A)| \text{vol})$. However, in terms of the description of $\overline{K}^0_R(\text{pt.})$ in Proposition 5, we have

$$[(V, |\rho_*(A)| \text{vol})] = [(V, \text{vol})] + a(\ln |\rho_*(A)|). \tag{28}$$

It follows that $a \circ \ln |\rho_*| = 0$. Similarly, if $r \in \text{Ker}(a)$ then there is a projective right-$R$-module $V$ and an automorphism $A$ of $V$ such that $r = \ln |\rho_*(A)|$. ∎

**Example 1** : If $R = \mathbf{Z}$ and $\rho : \mathbf{Z} \to \text{End}(\mathbf{C})$ is the standard inclusion then $\overline{K}^0_R(\text{pt.}) = \mathbf{R} \oplus \mathbf{Z}$. Specifically, if $V$ is a finitely-generated abelian group, let $V_{tor}$ be the torsion subgroup and let $\{e_i\}_{i=1}^m$ be an integral basis for $V/V_{tor}$. Given a volume form vol on $V_{\mathbf{C}}$, the isomorphism takes $(V, \text{vol})$ to $(\ln |V_{tor}| - \ln \text{vol}(e_1 \otimes_\rho 1, \ldots, e_m \otimes_\rho 1), m)$. If $R = \mathbf{C}$ and $\rho : \mathbf{C} \to \text{End}(\mathbf{C})$ is the standard map then $\overline{K}^0_R(\text{pt.}) = \mathbf{Z}$. This is because complex vector spaces, equipped with volume forms, are classified up to isomorphism by rank.

## 2.2 Analytic Pushforward of $\overline{K}$-Groups

Let $Z \to M \xrightarrow{\pi} B$ be a smooth fiber bundle with connected base $B$ and connected closed fibers $Z_b = \pi^{-1}(b)$. Let $F$ be a local system of right-$R$-modules on $M$. Let $H(Z; F|_Z)$ denote the $\mathbf{Z}$-graded local system of right-$R$-modules on $B$ whose fiber over $b \in B$ is isomorphic to the cohomology group $\text{H}^*(Z_b, F|_{Z_b})$. Define $[H(Z; F|_Z)] \in K^0_R(B)$ by

$$[H(Z; F|_Z)] = \sum_p (-1)^p [H^p(Z; F|_Z)]. \tag{29}$$

If we have a short exact sequence (2) of local systems on $M$ then we obtain a long exact sequence of local systems on $B$

$$\ldots \to H^*(Z; F^1|_Z) \to H^*(Z; F^2|_Z) \to H^*(Z; F^3|_Z) \to H^{*+1}(Z; F^1|_Z) \to \ldots \tag{30}$$

We use the notation of Subsection A.4.

**Definition 7** *The pushforward in real cohomology, denoted $\pi_! : H^*(M; \mathbf{R}) \to H^*(B; \mathbf{R})$, is given by*

$$\pi_!(\tau) = \int_Z e(TZ) \cup \tau. \tag{31}$$

**Definition 8** *The pushforward in $K_R^0$, denoted $\pi_! : K_R^0(M) \to K_R^0(B)$, is generated by*

$$\pi_!([F]) = [H(Z; F|_Z)]. \tag{32}$$

It follows from (30) that $\pi_!$ is well-defined on $K_R^0(M)$.

Pick a horizontal distribution $T^H M$ and a vertical Riemannian metric $g^{TZ}$ on the fiber bundle. If $F$ is a local system on $M$ and $h^F$ is a Hermitian metric on $F_{\mathbf{C}}$, let $h^H$ denote the $L^2$-metric on $H(Z; F|_Z)_{\mathbf{C}}$.

**Definition 9** *The pushforward in $\widehat{K}_R^0$, denoted $\pi_! : \widehat{K}_R^0(M) \to \widehat{K}_R^0(B)$, is generated by*

$$\pi_!([(F, h^F, \eta)]) = \left[ \left( H(Z; F|_Z), h^H, \int_Z e\left(TZ, \nabla^{TZ}\right) \wedge \eta - \mathcal{T}(T^H M, g^{TZ}, h^F) \right) \right]. \tag{33}$$

**Proposition 7** *The pushforward in $\widehat{K}_R^0$ is well-defined.*

**Pf.** : Suppose that we have a $\widehat{K}_R^0(M)$-relation in the sense of Definition 3. Let $T_f(F)$ be the torsion form of (2) and let $T_f(H)$ be the torsion form of (30). For $j \in \{1, 2, 3\}$, put $\mathcal{T}^j = \mathcal{T}(T^H M, g^{TZ}, h^{F^j})$. Define $\sigma \in \Omega^{even}(B)/\text{im}(d)$ by

$$\sigma = \mathcal{T}^2 - \mathcal{T}^1 - \mathcal{T}^3 + T_f(H) - \int_Z e\left(TZ, \nabla^{TZ}\right) \wedge T_f(F). \tag{34}$$

Then we must show that $\sigma = 0$.

**Lemma 1** *$\sigma$ is independent of $\{h^{F^j}\}_{j=1}^3$.*

**Pf.** : Let $\{h^{F^j}\}_{j=1}^3$ and $\{h'^{F^j}\}_{j=1}^3$ be two choices of Hermitian metrics on $\{F_\mathbf{C}^j\}_{j=1}^3$. Then there is a smooth family $\{h_s^{F^j}\}_{j=1}^3$ of Hermitian metrics, parametrized by $s \in [0,1]$, such that $h_0^{F^j} = h^{F^j}$ and $h_1^{F^j} = h'^{F^j}$. Put $\widetilde{M} = [0,1] \times M$ and $\widetilde{B} = [0,1] \times B$. Let $\pi_M : \widetilde{M} \to M$ and $\pi_B : \widetilde{B} \to B$ be the projections onto the second factors. Let $\widetilde{\pi} : \widetilde{M} \to \widetilde{B}$ be the projection $Id_{[0,1]} \times \pi$, with fiber $\widetilde{Z} \cong Z$. There is a natural horizontal distribution $T^H \widetilde{M} = \mathbf{R} \times T^H M$ and vertical Riemannian metric $g^{T\widetilde{Z}} = \pi_M^* g^{TZ}$ on $\widetilde{M}$. Put $\widetilde{F}^j = \pi_M^* F^j$. One has an equality of pairs

$$\left( H(\widetilde{Z}; \widetilde{F}^j|_{\widetilde{Z}}), \nabla^{H(\widetilde{Z};\widetilde{F}^j|_{\widetilde{z}})} \right) = \pi_B^* \left( H(Z; F^j|_Z), \nabla^{H(Z;F^j|_Z)} \right). \qquad (35)$$

We will abbreviate the left-hand-side of (35) by $(\widetilde{H}^j, \nabla^{\widetilde{H}^j})$. There is a Hermitian metric $h^{\widetilde{F}^j}$ on $\widetilde{F}^j$ whose restriction to $\{s\} \times M$ is $h_s^{F^j}$.

The exterior differentiation on $\widetilde{B}$ is given by

$$\widetilde{d} = ds \wedge \partial_s + d. \qquad (36)$$

Consider the torsion form $\mathcal{T}(T^H \widetilde{M}, g^{T\widetilde{Z}}, h^{\widetilde{F}^j}) \in \Omega^{even}(\widetilde{B})$. By Proposition 45,

$$\widetilde{d}\mathcal{T}(T^H \widetilde{M}, g^{T\widetilde{Z}}, h^{\widetilde{F}^j}) = \int_{\widetilde{Z}} e(T\widetilde{Z}, \nabla^{T\widetilde{Z}}) \wedge f(\nabla^{\widetilde{F}^j}, h^{\widetilde{F}^j}) - f(\nabla^{\widetilde{H}^j}, h^{\widetilde{H}^j}). \qquad (37)$$

By construction,

$$e(T\widetilde{Z}, \nabla^{T\widetilde{Z}}) = \pi_M^* e(TZ, \nabla^{TZ}). \qquad (38)$$

Equations (36) and (37) give that modulo $\mathrm{im}(d)$,

$$\frac{\partial \mathcal{T}^j}{\partial s} = \int_Z e(TZ, \nabla^{TZ}) \wedge i_{\partial_s} f(\nabla^{\widetilde{F}^j}, h^{\widetilde{F}^j}) - i_{\partial_s} f(\nabla^{\widetilde{H}^j}, h^{\widetilde{H}^j}). \qquad (39)$$

On the other hand, the same type of arguments applied to Proposition 41 give that modulo $\mathrm{im}(d)$,

$$\frac{\partial T_f(F)}{\partial s} = \sum_{j=1}^3 (-1)^j \, i_{\partial_s} f(\nabla^{\widetilde{F}^j}, h^{\widetilde{F}^j}) \qquad (40)$$

and

$$\frac{\partial T_f(H)}{\partial s} = \sum_{j=1}^3 (-1)^j \, i_{\partial_s} f(\nabla^{\widetilde{H}^j}, h^{\widetilde{H}^j}). \qquad (41)$$

Combining equations (39), (40) and (41), we see that $\sigma$ is independent of $s$.
∎

We now continue with the proof of Proposition 7. As $\sigma$ is defined solely in terms of complex vector bundles, we may as well assume that $R = \mathbf{C}$ and $\rho$ is the identity map. If the short exact sequence (2) splits, as a sequence of flat vector bundles, then it is easy to see that $\sigma$ vanishes. Choose a splitting $\mu : F^3 \to F^2$ of (2), as a sequence of smooth topological vector bundles.

Using $\mu$, there is an isomorphism of smooth topological vector bundles

$$F^2 \cong F^1 \oplus F^3. \tag{42}$$

In terms of this decomposition, we can write

$$\nabla^{F^2} = \begin{pmatrix} \nabla^{F^1} & \alpha \\ 0 & \nabla^{F^3} \end{pmatrix}, \tag{43}$$

with $\alpha \in \Omega^1(M; \mathrm{Hom}(F^3, F^1))$. The flatness of $\nabla^{F^2}$ is equivalent to the flatness of $\nabla^{F^1}$ and $\nabla^{F^3}$, along with $\alpha$ being covariantly-constant, i.e. in supernotation,

$$\nabla^{F^1}\alpha + \alpha\nabla^{F^3} = 0. \tag{44}$$

Following [9, Appendix A], we now show how to effectively rescale $\alpha$.

By Lemma 1, we may use arbitrary Hermitian metrics on the $\{F^j\}_{j=1}^3$. Choose Hermitian metrics $h^{F^1}$ and $h^{F^3}$ and take $h^{F^2} = h^{F^1} \oplus h^{F^3}$ so that (42) becomes an isometry. For $\epsilon \in (0,1]$, put

$$h_\epsilon^{F^1} = h^{F^1}, \quad h_\epsilon^{F^2} = h^{F^1} \oplus \epsilon^{-1} h^{F^3}, \quad h_\epsilon^{F^3} = \epsilon^{-1} h^{F^3}. \tag{45}$$

Let $X^1$ and $X^3$ be the $X$'s of (266) for the flat bundles $F^1$ and $F^3$, defined using $h^{F^1}$ and $h^{F^3}$, and define $\alpha^* \in \Omega^1(M; \mathrm{Hom}(F^1, F^3))$, the adjoint to $\alpha$, using $h^{F^1}$ and $h^{F^3}$. Then the $X$ of (266) for the $\mathbf{Z}$-graded vector bundle $F^1 \oplus F^2 \oplus F^3$, defined using $\{h_\epsilon^{F^j}\}_{j=1}^3$, takes the form

$$X_\epsilon = \begin{pmatrix} X^1 & 0 & 0 & 0 \\ 0 & X^1 & -\alpha/2 & 0 \\ 0 & \epsilon\alpha^*/2 & X^3 & 0 \\ 0 & 0 & 0 & X^3 \end{pmatrix}. \tag{46}$$

Putting
$$q = \begin{pmatrix} 0 & 0 & 0 & 0 \\ 0 & 0 & 0 & 0 \\ 0 & 0 & I & 0 \\ 0 & 0 & 0 & I \end{pmatrix}, \tag{47}$$

one has
$$\epsilon^{-q/2} X_\epsilon \, \epsilon^{q/2} = \begin{pmatrix} X^1 & 0 & 0 & 0 \\ 0 & X^1 & -\sqrt{\epsilon}\,\alpha/2 & 0 \\ 0 & \sqrt{\epsilon}\,\alpha^*/2 & X^3 & 0 \\ 0 & 0 & 0 & X^3 \end{pmatrix}. \tag{48}$$

Thus the torsion form $T_f(A', h_\epsilon^F)$ of the sequence (2) equals the torsion form $T_f(A', h^F)$ defined using the flat connection
$$\nabla_\epsilon^{F^2} = \begin{pmatrix} \nabla^{F^1} & \sqrt{\epsilon}\,\alpha \\ 0 & \nabla^{F^3} \end{pmatrix}. \tag{49}$$

A similar argument applies to the analytic torsion forms $\{\mathcal{T}^j\}_{j=1}^3$ and to the cohomology sequence (30). Using Lemma 1, the conclusion is that $\sigma$ is independent of the choice of $\epsilon \in (0,1]$ in (49). As in [9, Appendix A], one can show that $\sigma$ extends continuously to a function of $\epsilon \in [0,1]$. But when $\epsilon = 0$ then one is in the split situation and so $\sigma$ vanishes identically. ∎

The pushforward $\pi_! : \widehat{K}_R^0(M) \to \widehat{K}_R^0(B)$ depends explicitly on $T^H M$ and $g^{TZ}$ and hence has no topological meaning. We now show that when restricted to $\overline{K}_R^0(M)$, the pushforward is topological in nature.

**Proposition 8** *The pushforward in $\widehat{K}_R^0$ restricts to a pushforward*
$$\pi_! : \overline{K}_R^0(M) \to \overline{K}_R^0(B). \tag{50}$$

**Pf.** : Define $\pi_! : \Omega(M) \to \Omega(B)$ by
$$\pi_!(\eta) = \int_Z e\left(TZ, \nabla^{TZ}\right) \wedge \eta. \tag{51}$$

It is enough to show that there is a commutative diagram
$$\begin{array}{ccc} \widehat{K}_R^0(M) & \stackrel{c'}{\to} & \Omega^{odd}(M) \\ \pi_! \downarrow & & \pi_! \downarrow \\ \widehat{K}_R^0(B) & \stackrel{c'}{\to} & \Omega^{odd}(B). \end{array} \tag{52}$$

This follows from Proposition 45. ∎

**Proposition 9** *The pushforward in $\overline{K}_R^0$ is independent of the choices of $T^H M$ and $g^{TZ}$.*

**Pf.** : Let $(T^H M, g^{TZ})$ and $(T'^H M, g'^{TZ})$ be two choices of horizontal distributions and vertical Riemannian metrics on the fiber bundle. Then there is a smooth family $(T_s^H M, g_s^{TZ})$ of pairs, parametrized by $s \in [0,1]$, such that $(T_0^H M, g_0^{TZ}) = (T^H M, g^{TZ})$ and $(T_1^H M, g_1^{TZ}) = (T'^H M, g'^{TZ})$. Define $\widetilde{M}$ and $\widetilde{B}$ as in the proof of Lemma 1. Define $i_0 : B \to \widetilde{B}$ by $i_0(b) = (0,b)$ and $i_1 : B \to \widetilde{B}$ by $i_1(b) = (1,b)$. There is a horizontal distribution $T^H \widetilde{M}$ whose restriction to $\{s\} \times M$ is $\mathbf{R} \times T_s^H M$ and a vertical Riemannian metric $g^{T\widetilde{Z}}$ whose restriction to $\{s\} \times M$ is $g_s^{TZ}$.

Given a finite set of generators $\mathcal{F}^j = (F^j, h^{F^j}, \eta^j)$ in $\widehat{K}_R^0(M)$, put $\widetilde{\mathcal{F}}^j = \pi_M^* \mathcal{F}^j$. If $\sum_j n_j [\mathcal{F}^j]$ lies in $\overline{K}_R^0(M)$ then $\sum_j n_j \widetilde{\mathcal{F}}^j$ lies in $\overline{K}_R^0(\widetilde{M})$ and

$$\widetilde{k} = \widetilde{\pi}_! \sum_j n_j \widetilde{\mathcal{F}}^j \tag{53}$$

lies in $\overline{K}_R^0(\widetilde{B})$. By construction, $i_0^* \widetilde{k}$ is the pushforward of $\sum_j n_j [\mathcal{F}^j]$ using $(T^H M, g^{TZ})$ and $i_1^* \widetilde{k}$ is the pushforward of $\sum_j n_j [\mathcal{F}^j]$ using $(T'^H M, g'^{TZ})$. The proposition now follows from Proposition 4. ∎

**Proposition 10** *Let $R$ and $R'$ be two rings obeying our assumptions. Let $Z \to M \xrightarrow{\pi} B$ and $Z' \to M' \xrightarrow{\pi'} B$ be two fiber bundles over $B$ with connected closed fibers. Let $Z'' \to M'' \xrightarrow{\pi''} B$ be the product fiber bundle, with fiber $Z'' = Z \times Z'$. Let $p : M'' \to M$ and $p' : M'' \to M'$ be the natural projection maps. Then for $k \in \overline{K}_R^0(M)$ and $k' \in \overline{K}_{R'}^0(M')$, one has an identity in $\overline{K}_{R \otimes R'}^0(B)$*

$$\pi''_!(p^* k \cdot p'^* k') = (\pi_! k) \cdot (\pi'_! k'). \tag{54}$$

**Pf.** : Using the product formula for the analytic torsion forms given in [9, Prop. 3.28], this follows from a straightforward computation. We omit the details. ∎

**Remark 1 :** An important special case of the pushforward arises when there is an element $k \in \overline{K}_R^0(M)$ and an explicit trivialization of $b(\pi_! k)$ in $K_R^0(B)$. Then there is a canonical lifting of $\pi_! k \in \overline{K}_R^0(B)$ to $\mathrm{H}^{even}(B; \mathbf{R})$ in the sequence
$$\mathrm{H}^{even}(B; \mathbf{R}) \xrightarrow{a} \overline{K}_R^0(B) \xrightarrow{b} K_R^0(B). \tag{55}$$
In particular, suppose that $F$ is a local system of right-$R$-modules on $M$, $F_{\mathbf{C}}$ admits a covariantly-constant Hermitian metric $h^F$ and $H(Z; F|_Z) = 0$. Then we can take $\eta = 0$, $k = [(F, h^F, 0)]$ and pullback $\pi_! k$ to $-\mathcal{T} \in \mathrm{H}^{even}(B; \mathbf{R})$. Presumably this is the same up to constants as the higher Reidemeister torsion defined under the same circumstances by Igusa and Klein [19, 21].

We now look at the case when $B$ is a point.

**Proposition 11** *If $Z$ is a connected closed manifold, consider the fiber bundle $Z \xrightarrow{\pi} \mathrm{pt}$. Let $z_0$ be a point in $Z$ and let $s : \mathrm{pt.} \to Z$ be the section of $\pi$ given by $s(\mathrm{pt.}) = z_0$. Then $\pi_! : \overline{K}_R^0(Z) \to \overline{K}_R^0(\mathrm{pt.})$ is given by*
$$\pi_! = \chi(Z)\, s^*. \tag{56}$$

**Pf. :** Let $\sum_j n_j [(F^j, h^{F^j}, \eta^j)]$ be an element of $\overline{K}_R^0(Z)$. For each $j$, put
$$H(j) = H(Z; F^j\big|_Z). \tag{57}$$
Note that for each $j$, $H(j)$ is $\mathbf{Z}$-graded. In what follows, many of the objects will be implicitly graded. Choose a Riemannian metric $g^{TZ}$ on $Z$. Letting $\eta_{[0]}^j$ denote the 0-form part of $\eta^j$, we have
$$\pi_! \sum_j n_j [(F^j, h^{F^j}, \eta^j)] = \sum_j n_j [(H(j), h^{H(j)}, \int_Z e(\nabla^{TZ}) \wedge \eta_{[0]}^j - \mathcal{T}^j)]. \tag{58}$$

Put $\Gamma = \pi_1(Z, z_0)$. If $V^j$ is the fiber of $F^j$ over $z_0$, we can consider $V^j$ to be a left-$\mathbf{Z}\Gamma$-module and a right-$R$-module. Take a CW-decomposition of $Z$. Then the space of $F^j$-valued cellular cochains on $Z$ is
$$C^*(Z; F^j) = \mathrm{Hom}_{\mathbf{Z}\Gamma}(C_*(\widetilde{Z}), V^j), \tag{59}$$
where $C_*(\widetilde{Z})$ is a free left-$\mathbf{Z}\Gamma$-module with a basis given by a lifting of the cells of $Z$. Let $C$ be the $R$-cochain complex
$$C = \bigoplus_j n_j\, C^*(Z; F^j). \tag{60}$$

Put $C_{\mathbf{C}} = C \otimes_\rho \mathbf{C}^n$. As $\sum_j n_j [(F^j, h^{F^j}, \eta^j)]$ lies in $\overline{K}_R^0(Z)$, we have

$$\sum_j n_j \, d \, \log(\mathrm{vol}(h^{F^j})) = \sum_j n_j \, c_1(\nabla^{F^j}, h^{F^j}) = d \sum_j n_j \, \eta^j_{[0]}, \tag{61}$$

which implies that the complex line bundle $\otimes_j \left(\Lambda^{max} F_{\mathbf{C}}^j\right)^{n_j}$ on $Z$ has a covariantly-constant volume form $\otimes_j \left(e^{-\eta^j_{[0]}} \mathrm{vol}(h^{F^j})\right)^{n_j}$. Equivalently, the action of $\mathbf{Z}\Gamma$ on the complex line $L = \otimes_j \left(\Lambda^{max} V_{\mathbf{C}}^j\right)^{n_j}$ preserves a certain volume form. Using the cellular basis of $C^*(\widetilde{Z})$ and the volume form on $L$, $C_{\mathbf{C}}$ acquires a volume form $\mathrm{vol}(C)$.

Define a right-$R$-module $H$ by

$$H = \bigoplus_j n_j \, H(j). \tag{62}$$

Letting $H(j)_{\mathbf{C}}$ have the volume form

$$\mathrm{vol}(H(j)) = e^{-\int_Z e(TZ, \nabla^{TZ}) \wedge \eta^j_{[0]}} \mathrm{vol}(h^{H(j)}), \tag{63}$$

$H_{\mathbf{C}}$ acquires a volume form $\mathrm{vol}(H)$. Now $\mathcal{T}^j \in \mathbf{R}$ is given by the Ray-Singer analytic torsion [31]

$$\mathcal{T} = -\frac{1}{2} \frac{d}{ds}\bigg|_{s=0} \sum_{p=0}^{dim(Z)} (-1)^p \, p \, \mathrm{Tr}\left(\triangle'_p\right)^{-s}, \tag{64}$$

By the Cheeger-Müller theorem [13, 28, 29], $\sum_j n_j \mathcal{T}^j$ is the same as the Reidemeister torsion $T$ of $C$, computed using $\mathrm{vol}(C)$ and $\mathrm{vol}(H)$. Thus in the notation of Propositions 5 and 6,

$$\pi_! \sum_j n_j [(F^j, h^{F^j}, \eta^j)] = -a(T) + \sum_j n_j [(H(j), \mathrm{vol}(H(j)))]. \tag{65}$$

In fact, the right-hand-side of (65) is independent of the choice of volume forms $\{\mathrm{vol}(H(j))\}$, provided that the same volume forms are used to compute $T$. The result now follows from the following algebraic proposition. ∎

**Proposition 12** *Let $X$ be a finite chain complex of based finitely-generated free left-$\mathbf{Z}\Gamma$-modules. Let $V$ be an $R\Gamma$-module which is finitely-generated over $R$. Define cochain complexes by $C = \mathrm{Hom}_{\mathbf{Z}\Gamma}(X, V)$ and $C_{\mathbf{C}} = \mathrm{Hom}_{\mathbf{Z}\Gamma}(X, V_{\mathbf{C}})$. Let $H$ be the cohomology of $C$ and let $H_{\mathbf{C}}$ be the cohomology of $C_{\mathbf{C}}$. Let $\mathrm{vol}(V)$ be a $\Gamma$-invariant volume form on $V_{\mathbf{C}}$. Let $\mathrm{vol}(C)$ be the volume form on $C_{\mathbf{C}}$ constructed from $\mathrm{vol}(V)$ and the given basis of $X$. Choose volume forms $\{\mathrm{vol}(H^p)\}$ on $\{H_{\mathbf{C}}^p\}$. Let $T$ be the Reidemeister torsion of $(C, H)$, defined using $\mathrm{vol}(C)$ and $\mathrm{vol}(H)$. Then in $\overline{K}_R^0(\mathrm{pt.})$, one has*

$$- a(T) + \sum_p (-1)^p [(H^p, \mathrm{vol}(H^p))] = \chi(X) \cdot [(V, \mathrm{vol}(V))]. \qquad (66)$$

**Pf. :** (due to Wolfgang Lück) We use the following lemma, which is essentially a result of Milnor [27] (see also [25]).

**Lemma 2** *Let $0 \to B_{\mathbf{C}} \to C_{\mathbf{C}} \to D_{\mathbf{C}} \to 0$ be a short exact sequence of finite $\mathbf{C}$-cochain complexes. Let $H_{\mathbf{C}}$ be the acyclic complex of the long cohomology sequence. Suppose that the cochain groups and cohomology groups of $B_{\mathbf{C}}$, $C_{\mathbf{C}}$ and $D_{\mathbf{C}}$ are equipped with volume forms. Let $T(B)$, $T(C)$, $T(D)$ and $T(H)$ be the corresponding Reidemeister torsions. Let $T(B^p, C^p, D^p)$ be the Reidemeister torsion of the short exact sequence*

$$0 \to B_{\mathbf{C}}^p \to C_{\mathbf{C}}^p \to D_{\mathbf{C}}^p \to 0. \qquad (67)$$

*Then*

$$T(B) - T(C) + T(D) = T(H) + \sum_p (-1)^p \, T(B^p, C^p, D^p). \qquad (68)$$

**Lemma 3** *Let $C$ be a finite acyclic right-$R$-cochain complex. Let $\{\mathrm{vol}(C^p)\}$ be volume forms on $\{C_{\mathbf{C}}^p\}$. Let $T$ be the Reidemeister torsion of $C$. Then in $\overline{K}_R^0(\mathrm{pt.})$, one has*

$$- a(T) + \sum_p (-1)^p [(C^p, \mathrm{vol}(C^p))] = 0. \qquad (69)$$

**Pf. :** We do induction on the length $l$ of $C$. If $l \le 2$ then $C$ is a short exact sequence and the lemma is easy to check. Suppose that $l \ge 3$ and that the lemma is true for $l - 1$. Let $c_*$ be the coboundary operators of $C$. Let $D$ be the cochain complex concentrated in degrees $l - 1$ and $l$, and given there by the identity map $C^l \to C^l$. Put $\text{vol}(D^{l-1}) = \text{vol}(D^l) = \text{vol}(C^l)$. There is a cochain map $r : C \to D$ so that $r_{l-1} = c_{l-1}$ and $r_l = Id$. If $B = \ker(r)$, we have a short exact sequence of acyclic R-cochain complexes $0 \to B \to C \to D \to 0$. Equip $\{B_{\mathbf{C}}^p\}$ with volume forms. From Lemma 2,

$$T(C) = T(B) - \sum_p (-1)^p\, T(B^p, C^p, D^p). \tag{70}$$

Applying the induction hypothesis to the acyclic cochain complexes $B$ and $0 \to B^p \to C^p \to D^p \to 0$, we obtain

$$\begin{aligned}
a(T(C)) &= a(T(B)) - \sum_p (-1)^p\, a(T(B^p, C^p, D^p)) \\
&= \sum_p (-1)^p\, [(B^p, \text{vol}(B^p))] \\
&\quad - \sum_p (-1)^p\, ([(B^p, \text{vol}(B^p))] - [(C^p, \text{vol}(C^p))] + [(D^p, \text{vol}(D^p))]) \\
&= \sum_p (-1)^p\, [(C^p, \text{vol}(C^p))]. \tag{71}
\end{aligned}$$

This proves the lemma. ∎

Let $\sigma(X)$ denote the left-hand-side of (66).

**Lemma 4** *Let $0 \to X_1 \to X_2 \to X_3 \to 0$ be a based short exact sequence of finite based free left-$\mathbf{Z}\Gamma$-chain complexes. Let $V$ be an $R\Gamma$-module which is finitely-generated over $R$. Let $\text{vol}(V)$ be a $\Gamma$-invariant volume form on $V_{\mathbf{C}}$. For $i \in \{1, 2, 3\}$, put $C_i = \text{Hom}_{\mathbf{Z}\Gamma}(X_i, V)$. Let $\text{vol}(C_i)$ be the volume form on $C_{i\mathbf{C}}$ constructed from $\text{vol}(V)$ and the given basis of $X_i$. Then for any choice of volume forms $\{\text{vol}(H^*(C_i))\}_{1 \le i \le 3}$, in $\overline{K}_R^0(\text{pt.})$ one has*

$$\sigma(X_1) - \sigma(X_2) + \sigma(X_3) = 0. \tag{72}$$

**Pf.** : We have a volume-preserving short exact sequence of **C**-cochain complexes

$$0 \to C_{1\mathbf{C}} \to C_{2\mathbf{C}} \to C_{3\mathbf{C}} \to 0. \tag{73}$$

Let $H_{\mathbf{C}}$ denote its long cohomology sequence. From Lemma 2,

$$T(C_{1\mathbf{C}}) - T(C_{2\mathbf{C}}) + T(C_{3\mathbf{C}}) = T(H). \tag{74}$$

From Lemma 3,

$$\begin{aligned} a(T(H)) &= \sum_p (-1)^p \left[(H^p(C_1), \mathrm{vol}(H^p(C_1))\right] \\ &\quad - \sum_p (-1)^p \left[(H^p(C_2), \mathrm{vol}(H^p(C_2))\right] \\ &\quad + \sum_p (-1)^p \left[(H^p(C_3), \mathrm{vol}(H^p(C_3))\right]. \end{aligned} \qquad (75)$$

The lemma follows from combining (74) and (75). ∎

We now prove the proposition by induction on the length $l$ of $X$. The case $l = 0$ is easy to check. Suppose that $l \geq 1$ and that the proposition is true for $l-1$. Let $X|_{l-1}$ be the initial segment of $X$ of length $l-1$. Let $X(l)$ be the cochain complex given by just the final term in $X$. Then there is a based short exact sequence of finite based free left-$\mathbf{Z}\Gamma$-chain complexes

$$0 \to X(l) \to X \to X|_{l-1} \to 0. \qquad (76)$$

By Lemma 4, for any choice of volume forms on the cohomology groups,

$$\sigma(X) = \sigma(X(l)) + \sigma(X|_{l-1}). \qquad (77)$$

Applying the induction hypothesis to $X(l)$ and $X|_{l-1}$, the proposition follows. ∎

**Example 2 :** Suppose that $R = \mathbf{Z}$ and $\rho : \mathbf{Z} \to \mathrm{End}(\mathbf{C})$ is the standard inclusion. Let $Z$ be a closed connected manifold. Let $F$ be the trivial local system on $Z$ with fiber $V = \mathbf{Z}$. Let $h^F$ be the flat Hermitian metric on $F_{\mathbf{C}}$ normalized to be 1 on the generators of $F \subset F_{\mathbf{C}}$. Choose a Riemannian metric $g^{TZ}$ on $Z$. Let $\mathrm{vol}_{L^2}(H^p)$ be the volume form on $H^p(Z; \mathbf{R})$ coming from the $L^2$-inner product. Let $\mathrm{vol}_{\mathbf{Z}}(H^p)$ be the volume form on $H^p(Z; \mathbf{R})$ coming from an integral basis of $H^p(Z; \mathbf{Z})/H^p(Z; \mathbf{Z})_{tor}$. Using the results of Example 1, the statement of Proposition 11 boils down to [13, Theorem 8.35], namely that the analytic torsion of $(Z, F)$ is given by

$$\mathcal{T} = \sum_p (-1)^p \left( \ln |H^p(Z; \mathbf{Z})_{tor}| - \ln \frac{\mathrm{vol}_{L^2}(H^p)}{\mathrm{vol}_{\mathbf{Z}}(H^p)} \right). \qquad (78)$$

**Remark 2 :** In the case of $Z \xrightarrow{\pi}$ pt., we have seen that only the degree-0 part of $\eta$ is relevant. More generally, if $\dim(B) \leq 2r$, one can redo the definitions truncating $\eta$ at degree $2r$ and $c(\nabla^F, h^F)$ at degree $2r+1$.

## 2.3 Possible Topological Interpretation

With the notation of Subsection 2.1, we have a commutative diagram

$$
\begin{array}{ccccccc}
\mathrm{H}^{even}(M;\mathbf{R}) & \xrightarrow{a} & \overline{K}^0_R(M) & \xrightarrow{b} & K^0_R(M) & \xrightarrow{c} & \mathrm{H}^{odd}(M;\mathbf{R}) \\
\pi_! \downarrow & & \pi_! \downarrow & & \pi_! \downarrow & & \pi_! \downarrow \\
\mathrm{H}^{even}(B;\mathbf{R}) & \xrightarrow{a} & \overline{K}^0_R(B) & \xrightarrow{b} & K^0_R(B) & \xrightarrow{c} & \mathrm{H}^{odd}(B;\mathbf{R}).
\end{array} \tag{79}
$$

Let $K_{alg}(R) = K_0(R) \times BGL(R)^+_\delta$ be the classifying space for algebraic K-theory, where $\delta$ denotes the discrete topology on $GL(R)$ and $+$ denotes Quillen's plus construction. Let $K^*_{alg,R}$ be the corresponding generalized cohomology theory, so that $K^0_{alg,R}(M) = [M, K_{alg}(R)]$. Consider the map

$$K_{alg}(R) \xrightarrow{\rho_*} K_{alg}(\mathbf{C}) \xrightarrow{\beta} \prod_{j=1}^{\infty} K(\mathbf{R}, 2j+1), \tag{80}$$

where $\beta$ is the map given by the Borel regulator classes. Let $\mathcal{F}_R$ be the homotopy fiber of $\beta \circ \rho_*$, a map of infinite loop spaces, and let $\mathcal{F}^*_R$ be the corresponding generalized cohomology theory. The homotopy exact sequence gives

$$\ldots \to K_1(R) \to \mathbf{R} \to \mathcal{F}^0_R(\mathrm{pt.}) \to K_0(R) \to 0, \tag{81}$$

which can be compared with Proposition 6.

Given a generalized cohomology theory $E$ and a fiber bundle $Z \to M \xrightarrow{\pi} B$ with connected compact fibers, the Becker-Gottlieb-Dold transfer gives a map $\mathrm{tr} : E^*(M) \to E^*(B)$ [4, 15]. There is a corresponding commutative diagram

$$
\begin{array}{ccccccc}
\mathrm{H}^{even}(M;\mathbf{R}) & \longrightarrow & \mathcal{F}^0_R(M) & \longrightarrow & K^0_{alg,R}(M) & \longrightarrow & \mathrm{H}^{odd}(M;\mathbf{R}) \\
\mathrm{tr} \downarrow & & \mathrm{tr} \downarrow & & \mathrm{tr} \downarrow & & \mathrm{tr} \downarrow \\
\mathrm{H}^{even}(B;\mathbf{R}) & \longrightarrow & \mathcal{F}^0_R(B) & \longrightarrow & K^0_{alg,R}(B) & \longrightarrow & \mathrm{H}^{odd}(B;\mathbf{R}).
\end{array} \tag{82}
$$

In the case of ordinary cohomology with real coefficients, tr is the same as the $\pi_!$ of Definition 7. Furthermore, with our regularity assumption on $R$, there is a map $K^0_R(M) \to K^0_{alg,R}(M)$ which essentially comes from the map $K_0(R) \times BGL(R)_\delta \to K_0(R) \times BGL(R)^+_\delta$. In terms of the pushforward $\pi_!$ of Definition 8, Dwyer and Williams have shown [16] that there is a commutative diagram

$$\begin{array}{ccc} K^0_R(M) & \to & K^0_{alg,R}(M) \\ \pi_! \downarrow & & \text{tr} \downarrow \\ K^0_R(B) & \to & K^0_{alg,R}(B). \end{array} \tag{83}$$

Superimposing (79) on (82), the natural guess is that the pushforward $\pi_! : \overline{K}^0_R(M) \to \overline{K}^0_R(B)$ is essentially the Becker-Gottlieb-Dold transfer.

**Conjecture 1** *There is a natural map $\overline{K}^0_R(M) \to \mathcal{F}^0_R(M)$ such that the following diagrams commute :*

$$\begin{array}{ccccc} \mathrm{H}^{even}(M;\mathbf{R}) & \to & \overline{K}^0_R(M) & \to & K^0_R(M) \\ \downarrow & & \downarrow & & \downarrow \\ \mathrm{H}^{even}(M;\mathbf{R}) & \to & \mathcal{F}^0_R(M) & \to & K^0_{alg,R}(M) \end{array} \tag{84}$$

*and*

$$\begin{array}{ccc} \overline{K}^0_R(M) & \to & \mathcal{F}^0_R(M) \\ \pi_! \downarrow & & \text{tr} \downarrow \\ \overline{K}^0_R(B) & \to & \mathcal{F}^0_R(B). \end{array} \tag{85}$$

Conjecture 1 is consistent with Proposition 11, as the Becker-Gottlieb-Dold transfer of the fiber bundle $Z \to \text{pt.}$ is given by $\chi(Z) \, s^*$.

One can compare Conjecture 1 with the corresponding result of [22], in which $\pi_!$ was an analytic index in $\mathbf{R}/\mathbf{Z}$ K-theory and tr was a topological index. The method of proof in [22] was to use the topological pairing

$$K_{-1}(B) \times K^{-1}(B; \mathbf{R}/\mathbf{Z}) \to \mathbf{R}/\mathbf{Z}, \tag{86}$$

which could be described analytically in terms of eta-invariants. This method of proof will not work here, because of the fact that the pushforward of Proposition 11 gives little information about $\overline{K}^0_R(Z)$.

**Remark 3 :** One can ask similar questions in the holomorphic case, following the work of Gillet and Soulé. Namely, if $M$ is an arithmetic variety, define $\widehat{K}_0(M)$ as in [17, Definition 6.1]. Let $\overline{K}_0(M)$ be the kernel of the map ch : $\widehat{K}_0(M) \to A(M)$. Suppose that $M \stackrel{\pi}{\to} B$ is a smooth projective map of arithmetic varieties such that $M_\infty \stackrel{\pi_\infty}{\to} B_\infty$ is a Kähler fibration. There is a pushforward $\widehat{K}_0(M) \to \widehat{K}_0(B)$ [18] which restricts to a pushforward $\pi_! : \overline{K}_0(M) \to \overline{K}_0(B)$. From [8, Theorem 3.10], $\pi_!$ is independent of the choice of the form $\omega$ on $M_\infty$ defining the Kähler fibration. (Strictly speaking, one should also relate the holomorphic torsion forms of [18] and [8].) One can then ask for a purely holomorphic description of $\pi_!$.

A related question occurs when $B$ is a point. Let $Z$ be a closed Hermitian complex manifold. Let $E_1$ and $E_2$ be two flat Hermitian vector bundles on $Z$ with $\text{rk}(E_1) = \text{rk}(E_2)$ such that for some $p$, $\text{H}^{p,*}(Z; E_1) = \text{H}^{p,*}(Z; E_2) = 0$. Then by [32], the difference of the holomorphic torsions $\mathcal{T}_p(Z, E_1) - \mathcal{T}_p(Z, E_2)$ is independent of the Hermitian metric on $Z$. A problem which is implicit in [32] is to give a purely holomorphic description of $\mathcal{T}_p(Z, E_1) - \mathcal{T}_p(Z, E_2)$. See [14], [32] and [33] for results along this line.

## 3  Flat Duality Bundles

We consider flat real vector bundles $E$ on a manifold $M$ which are equipped with a parallel nondegenerate $\epsilon$-symmetric bilinear form, where $\epsilon = \pm 1$. Given a reduction $J^E$ of the structure group of $E$ (as a topological vector bundle) to its maximal compact subgroup, in Subsection 3.1 we define a characteristic form $p(\nabla^E, J^E) \in \Omega^{4*+1-\epsilon}(M)$. In Subsection 3.2 we extend these considerations to the setting of superconnections and define the notion of a flat duality superconnection on $E$, along with its characteristic form $p(A', J^E) \in \Omega^{4*+1-\epsilon}(M)$. In Subsection 3.3 we define a group $\overline{L}^0_\epsilon(M)$ which essentially consists of pairs $(E, J^E)$ along with an explicit writing of $p(\nabla^E, J^E)$ as an exact form. We show that $\overline{L}^0_\epsilon(M)$ is a homotopy-invariant of $M$ and we show how to detect elements of $\overline{L}^0_\epsilon(M)$ by means of real-valued eta-invariants of tangential signature operators. In Subsection 3.4 we assume that $E$ is $\mathbf{Z}$-graded and use the ensuing rescaling of $J^E$ to construct the finite-dimensional eta-form $\widetilde{\eta}$ of [7]. In Subsection 3.5 we start with a smooth fiber bundle $M \stackrel{\pi}{\to} B$ with closed oriented fibers $Z$. Adding a horizontal distribution $T^H M$ and a vertical Riemannian metric $g^{TZ}$, we show that the infinite-dimensional vector bundle $W$ on $B$ of vertical differential forms acquires a flat duality superconnection. In Subsection 3.6 we construct an analytic pushforward $\pi_! : \overline{L}^0_\epsilon(M) \to \overline{L}^0_\epsilon(B)$. This pushforward involves the eta-form of the vertical signature operator. We show that $\pi_!$ is independent of $T^H M$ and $g^{TZ}$.

One could generalize the results of this section to the setting of modules over rings with anti-involutions, but for simplicity we only discuss vector spaces over the reals.

## 3.1 Characteristic Classes of Flat Duality Bundles

Let $V$ be a finite-dimensional real vector space. Let $(\cdot, \cdot)_V$ be a nondegenerate bilinear form on $V$. Given $\epsilon = \pm 1$, we say that $(\cdot, \cdot)_V$ is $\epsilon$-symmetric if $(v_1, v_2)_V = \epsilon(v_2, v_1)_V$ for all $v_1, v_2 \in V$. The automorphism group of $(V, (\cdot, \cdot)_V)$ is of the form $O(m, m')$ if $\epsilon = 1$ or $Sp(2m, \mathbf{R})$ if $\epsilon = -1$.

Let $M$ be a smooth connected manifold

**Definition 10** *A flat duality bundle on $M$ is a real vector bundle $E$ on $M$ with a flat connection $\nabla^E$ and a nondegenerate $\epsilon$-symmetric bilinear form $(\cdot, \cdot)_E$ as above which is covariantly-constant with respect to $\nabla^E$.*

Equivalently, we assume that the holonomy group of $E$ is a discrete subgroup of $O(m, m')$ (or $Sp(2m, \mathbf{R})$). We can choose an automorphism $J^E$ of $E$, considered as a smooth topological vector bundle, so that

- $(J^E)^2 = \epsilon$.

- $(v_1, v_2)_E = (J^E v_1, J^E v_2)_E$.

- The inner product $< v_1, v_2 >_E = (v_1, J^E v_2)_E$ is positive-definite.

The choice of $J^E$ reduces the structure group of $E$, as a smooth topological vector bundle, to $O(m) \times O(m')$ (or $U(m)$).

**Definition 11** *Define $\omega(\nabla^E, J^E) \in \Omega^1(M; \text{End}(E))$ by*

$$\omega(\nabla^E, J^E) = (J^E)^{-1}(\nabla^E J^E). \tag{87}$$

When no ambiguity can arise, we will abbreviate $\omega(\nabla^E, J^E)$ by $\omega$. One has

$$\nabla^E \omega = -\omega^2. \tag{88}$$

**Proposition 13** *The adjoint of $\nabla^E$ with respect to the inner product $< \cdot, \cdot >_E$ is*

$$\left(\nabla^E\right)^T = (J^E)^{-1}\nabla^E J^E = \nabla^E + \omega. \tag{89}$$

**Pf.** : Given sections $e_1$ and $e_2$ of $E$ and a vector field $X$ on $M$, we have

$$\begin{aligned}
< e_1, (\nabla^E)^T_X e_2 >_E &= X < e_1, e_2 >_E - < \nabla^E_X e_1, e_2 >_E \\
&= X(e_1, J^E e_2)_E - (\nabla^E_X e_1, J^E e_2)_E \\
&= (e_1, \nabla^E_X (J^E e_2))_E \\
&= < e_1, ((J^E)^{-1}\nabla^E_X J^E) e_2 >_E .
\end{aligned} \tag{90}$$

The proposition follows. ∎

**Definition 12** *Define a connection $\nabla^{E,u}$ on $E$ by*

$$\nabla^{E,u} = \nabla^E + \frac{\omega}{2}. \tag{91}$$

**Proposition 14** *The connection $\nabla^{E,u}$ preserves $< \cdot, \cdot >_E$ and commutes with $J^E$. Its curvature is given by*

$$\left(\nabla^{E,u}\right)^2 = -\frac{\omega^2}{4}. \tag{92}$$

**Pf.** : From Proposition 13, $\nabla^{E,u} = \frac{1}{2}\left(\nabla^E + (\nabla^E)^T\right)$ preserves $< \cdot, \cdot >_E$. Next,

$$\begin{aligned}
\nabla^{E,u} J^E &= \frac{1}{2}(\nabla^E + (J^E)^{-1}\nabla^E J^E)J^E = \frac{1}{2}(\nabla^E J^E + \epsilon\,(J^E)^{-1}\nabla^E) \\
&= \frac{1}{2}J^E((J^E)^{-1}\nabla^E J^E + \nabla^E) = J^E\,\nabla^{E,u}.
\end{aligned} \tag{93}$$

Finally,

$$\left(\nabla^{E,u}\right)^2 = \left(\nabla^E + \frac{1}{2}\omega\right)^2 = \frac{1}{2}\nabla^E\omega + \frac{\omega^2}{4} = -\frac{\omega^2}{4}. \tag{94}$$

The proposition follows. ∎

**Definition 13** *If $\epsilon = 1$, define $p(\nabla^E, J^E) \in \Omega^{4*}(M)$ by*

$$p(\nabla^E, J^E) = \text{tr}\left[J^E \cos\left(\frac{\omega^2}{8\pi}\right)\right]. \tag{95}$$

*If $\epsilon = -1$, define $p(\nabla^E, J^E) \in \Omega^{4*+2}(M)$ by*

$$p(\nabla^E, J^E) = -\text{tr}\left[J^E \sin\left(\frac{\omega^2}{8\pi}\right)\right]. \tag{96}$$

If $\epsilon = 1$, put

$$P_{\pm} = \frac{I \pm J^E}{2}, \quad E_{\pm} = P_{\pm} E, \quad \nabla^{E_{\pm}} = P_{\pm} \nabla^{E,u} P_{\pm}. \tag{97}$$

Then $\nabla^{E_{\pm}}$ is an orthogonal connection on the real vector bundle $E_{\pm}$.

If $\epsilon = -1$, put

$$P_{\pm} = \frac{I \mp i J^E}{2}, \quad E_{\pm} = P_{\pm}(E_{\mathbf{C}}), \quad \nabla^{E_{\pm}} = P_{\pm} \nabla^{E,u} P_{\pm}. \tag{98}$$

Then $\nabla^{E_{\pm}}$ is a Hermitian connection on the complex vector bundle $E_{\pm}$.

**Proposition 15** *We have*

$$p(\nabla^E, J^E) = \text{ch}(\nabla^{E_+}) - \text{ch}(\nabla^{E_-}). \tag{99}$$

**Pf.** : As $\nabla^{E,u}$ commutes with $P_{\pm}$, the curvature of $\nabla^{E_{\pm}}$ is given by

$$\left(\nabla^{E_{\pm}}\right)^2 = P_{\pm}\left(\nabla^{E,u}\right)^2 P_{\pm} = -\frac{1}{4} P_{\pm} \omega^2 P_{\pm}. \tag{100}$$

It follows that

$$\text{ch}(\nabla^{E_+}) - \text{ch}(\nabla^{E_-}) = \text{tr}\left[P_+ e^{\frac{\omega^2}{8i\pi}}\right] - \text{tr}\left[P_- e^{\frac{\omega^2}{8i\pi}}\right] \tag{101}$$

If $\epsilon = 1$ then we get

$$\text{ch}(\nabla^{E_+}) - \text{ch}(\nabla^{E_-}) = \text{tr}\left[J^E e^{\frac{\omega^2}{8i\pi}}\right] = \text{tr}\left[J^E \cos\left(\frac{\omega^2}{8\pi}\right)\right]. \tag{102}$$

If $\epsilon = -1$ then we get

$$\text{ch}(\nabla^{E_+}) - \text{ch}(\nabla^{E_-}) = \text{tr}\left[-i J^E e^{\frac{\omega^2}{8i\pi}}\right] = -\text{tr}\left[J^E \sin\left(\frac{\omega^2}{8\pi}\right)\right]. \tag{103}$$

The proposition follows. ∎

**Corollary 3** *The form $p(\nabla^E, J^E)$ is closed.*

**Pf. :** Using Proposition 15, this follows from the analogous property of the Chern character. ∎

We will show in Corollary 5 that the de Rham cohomology class $p(\nabla^E) \in \mathrm{H}^{4*+1-\epsilon}(M; \mathbf{R})$ of $p(\nabla^E, J^E)$ is independent of $J^E$. As in Subsection A.1, one can view $p(\nabla^E)$ as arising from the group cohomology of $O(m, m')$ (or $Sp(2m, \mathbf{R})$). Let $\mathrm{H}_c^*(O(m, m'); \mathbf{R})$ (or $\mathrm{H}_c^*(Sp(2m, \mathbf{R}); \mathbf{R})$) be the cohomology of the complex of Eilenberg-Maclane cochains on $O(m, m')$ (or $Sp(2m, \mathbf{R})$) which are continuous in their arguments. Then the inverse limit $\lim_{m,m' \to \infty} \mathrm{H}_c^*(O(m, m'); \mathbf{R})$ (or $\lim_{m \to \infty} \mathrm{H}_c^*(Sp(2m, \mathbf{R}); \mathbf{R})$) is a symmetric algebra with generators in degrees congruent to 0 (or 2) mod 4 [10]. There are forgetful maps

$$\mu : \mathrm{H}_c^*(O(m, m'); \mathbf{R}) \to \mathrm{H}^*(BO(m, m')_\delta; \mathbf{R}) \tag{104}$$

and

$$\mu : \mathrm{H}_c^*(Sp(2m, \mathbf{R}); \mathbf{R}) \to \mathrm{H}^*(BSp(2m, \mathbf{R})_\delta; \mathbf{R}). \tag{105}$$

The terms in $p(\nabla^E)$ are the pullbacks, under the classifying map of $E$ given by $\nu : M \to BO(m, m')_\delta$ (or $\nu : M \to BSp(2m, \mathbf{R})_\delta$), of the $\mu$-images of these generators.

### 3.2 Flat Duality Superconnections

Let $E$ be a real vector bundle on $M$ with an element $J^E \in \mathrm{Aut}(E)$ satisfying $(J^E)^2 = \epsilon$, $\epsilon = \pm 1$. Let $\nabla^{E,u}$ be a connection on $E$ which commutes with $J^E$. Given $S \in \Omega(M; \mathrm{End}(E))$, we can decompose $S$ as

$$S = \sum_{j \geq 0} S_j, \tag{106}$$

where $S_j$ is of partial degree $j$ in the Grassmann variables of $\Lambda(T^*M)$. We say that $S$ is even if for all $j \geq 0$,

$$J^E S_j = (-1)^j S_j J^E \tag{107}$$

and odd if for all $j \geq 0$,

$$J^E S_j = (-1)^{j+1} S_j J^E. \tag{108}$$

**Definition 14** *A $J^E$-superconnection on $E$ is an operator of the form $A = \nabla^{E,u} + S$ with $S$ odd.*

Clearly if we complexify and let $\tau = \frac{1}{\sqrt{\epsilon}} J^E$ define the $\mathbf{Z}_2$-grading on $E_{\mathbf{C}}$ then $A_{\mathbf{C}}$ is a superconnection on $E_{\mathbf{C}}$ in the sense of [30]. To get the signs right when manipulating $A$, one can use the sign conventions for the action of $A_{\mathbf{C}}$ on the $\mathbf{Z}_2$-graded vector space $E_{\mathbf{C}}$. By definition, the curvature of $A$ is $A^2$, a $C^\infty(M)$-linear endomorphism of $\Omega(M; E)$ which is given by multiplication by an even element of $\Omega(M; \operatorname{End}(E))$.

Now suppose that as in Subsection 3.1, $E$ is endowed with a nondegenerate $\epsilon$-symmetric bilinear form $(\cdot, \cdot)$. Extend $(\cdot, \cdot)_E$ to an $\Omega(M)$-valued bilinear form on $\Omega(M; E)$ by requiring that for $e, e' \in C^\infty(M; E)$ and $\omega \in \Omega^j(M)$, $\omega' \in \Omega^{j'}(M)$,

$$(\omega e, \omega' e')_E = (-1)^{\frac{j'(j'-1)}{2}} \omega \wedge \omega' \, (e, e')_E. \tag{109}$$

Let $\sigma$ be a new variable which commutes with $C^\infty(M; \operatorname{End}(E))$, anticommutes with $\Omega^1(M)$ and satisfies $\sigma^2 = 1$.

**Definition 15** *Given a $J^E$-superconnection $A$ and an even element $X$ of $\Omega(M; \operatorname{End}(E))$, we say that $(A, X)$ forms a flat pair if*

$$(A - X\sigma)^2 = 0. \tag{110}$$

*We say that $(A, X)$ is $(\cdot, \cdot)_E$-compatible if for all $e, e' \in C^\infty(M; E)$,*

$$\begin{aligned} d(e, e')_E &= (Ae, e')_E + (e, Ae')_E, \\ 0 &= (Xe, e')_E + (e, Xe')_E. \end{aligned} \tag{111}$$

*We say that $(A, X)$ is a flat duality superconnection if it is a $(\cdot, \cdot)_E$-compatible flat pair.*

Put

$$\begin{aligned} A' &= A - X, \\ A'' &= A + X. \end{aligned} \tag{112}$$

Here $A'$ is the sum of a connection $\nabla^E$ and an element of $\Omega(M; \operatorname{End}(E))$, and similarly for $A''$. One can show that (111) is equivalent to

$$d(e, e')_E = (A'e, e')_E + (e, A'e')_E. \tag{113}$$

Given $A'$ and $J^E$, we can recover $A$ and $X$ as the odd and even parts of $A'$, with respect to $J^E$. Note, however, that $A'_{\mathbf{C}}$ is inhomogeneous with respect to the $\mathbf{Z}_2$-grading on $E_{\mathbf{C}}$. For this reason, we generally write things in terms of $A$ and $X$.

In analogy with [9, Sections 1c,d], let $\cdot^T$ denote the linear map on $\Omega(M;\mathrm{End}(E))$ defined by the relations

- For $\alpha, \alpha' \in \Omega(M;\mathrm{End}(E))$,
$$(\alpha\alpha')^T = \alpha'^T \alpha^T. \tag{114}$$

- For $\omega \in \Omega^1(M)$,
$$\omega^T = -\omega. \tag{115}$$

- For $a \in C^\infty(M;\mathrm{End}(E))$, we have that $a^T$ is the transpose of $a$ in the ordinary sense, that is, for all $e, e' \in C^\infty(M;E)$,
$$<a^T e, e'>_E = <e, ae'>_E. \tag{116}$$

Given a $J^E$-superconnection $A = \nabla^{E,u} + S$, we define its transpose to be $A^T = (\nabla^{E,u})^T + S^T$. We say that a $J^E$-superconnection $A$ is symmetric if $A^T = A$.

**Proposition 16** *If the pair $(A, X)$ is $(\cdot, \cdot)_E$-compatible then $A$ is symmetric.*

**Pf. :** Suppose that $(A, X)$ is $(\cdot, \cdot)_E$-compatible. Let us write $A$ as $\nabla^{E,u} + S$, with $S_1 = 0$. From (111),
$$d(e, e')_E = (\nabla^{E,u} e, e')_E + (e, \nabla^{E,u} e')_E, \tag{117}$$

implying that for any vector field $V$ on $M$,
$$V<e, e'>_E = <\nabla^{E,u}_V e, e'>_E + <e, (J^E)^{-1} \nabla^{E,u}_V J^E e'>_E. \tag{118}$$

As $A$ is a $J^E$-superconnection, $(J^E)^{-1} \nabla^{E,u}_V J^E = \nabla^{E,u}_V$ and so (118) implies that $\nabla^{E,u}$ is symmetric. Furthermore, (111) implies that for $j \geq 0$,
$$0 = (S_j e, e')_E + (e, S_j e')_E. \tag{119}$$

Let $\{\omega_{j,r}\}$ be a local basis of $\Omega^j(M)$. Write
$$S_j = \sum_r \omega_{j,r} B_{j,r}, \tag{120}$$
with $B_{j,r} \in C^\infty(M; \text{End}(E))$ satisfying
$$J^E B_{j,r} = (-1)^{j+1} B_{j,r} J^E. \tag{121}$$
Then
$$\begin{aligned} 0 &= \sum_r [(\omega_{j,r} B_{j,r} e, e')_E + (e, \omega_{j,r} B_{j,r} e')_E] \\ &= \sum_r \omega_{j,r} \left[ (B_{j,r} e, e')_E + (-1)^{\frac{j(j-1)}{2}} (e, B_{j,r} e')_E \right]. \end{aligned} \tag{122}$$
Thus
$$0 = (B_{j,r} e, e')_E + (-1)^{\frac{j(j-1)}{2}} (e, B_{j,r} e')_E, \tag{123}$$
implying that
$$\begin{aligned} 0 &= <B_{j,r} e, e'>_E + (-1)^{\frac{j(j-1)}{2}} <e, (J^E)^{-1} B_{j,r} J^E e'>_E \\ &= <B_{j,r} e, e'>_E - (-1)^{\frac{j(j+1)}{2}} <e, B_{j,r} e'>_E. \end{aligned} \tag{124}$$
That is,
$$B_{j,r}^T = (-1)^{\frac{j(j+1)}{2}} B_{j,r}. \tag{125}$$
It follows that
$$\begin{aligned} S_j^T &= \sum_r (\omega_{j,r} B_{j,r})^T = \sum_r B_{j,r}^T \omega_{j,r}^T = (-1)^{\frac{j(j+1)}{2}} \sum_r B_{j,r}^T \omega_{j,r} \\ &= (-1)^{\frac{j(j+1)}{2}} \sum_r \omega_{j,r} B_{j,r}^T = \sum_r \omega_{j,r} B_{j,r} = S_j. \end{aligned} \tag{126}$$
Thus $A$ is symmetric. ∎

Let $(A, X)$ be a flat duality superconnection. Define $\phi$ as in (255).

**Definition 16** *Define $p(A', J^E) \in \Omega^{even}(M)$ by*
$$p(A', J^E) = \frac{1}{\sqrt{\epsilon}} \phi \operatorname{tr}\left[ J^E e^{-A^2} \right]. \tag{127}$$

**Proposition 17** *One has that $p(A', J^E)$ is a real closed element of $\Omega^{4*+1-\epsilon}(M)$.*

**Pf.** : Using the $\mathbf{Z}_2$-grading on $E_\mathbf{C}$ given by $\tau$, we have

$$p(A', J^E) = \text{ch}(A_\mathbf{C}). \tag{128}$$

It follows that $p(A', J^E)$ is real and closed. As $\text{tr}\left[J^E e^{-A^2}\right]$ is also real, $p(A', J^E)$ must be concentrated in degrees congruent to $(1 - \epsilon) \mod 4$. ∎

**Proposition 18** *Fix $(\cdot, \cdot)$ and $A'$. Let $s \in \mathbf{R}$ parametrize a smooth family of automorphisms $J^E(s)$ as above. Then*

$$\frac{\partial}{\partial s} p(A', J^E(s)) = d\left((2i\pi)^{-1/2} \frac{1}{2\sqrt{\epsilon}} \phi \, \text{tr}\left[J^E \left[(J^E)^{-1} \frac{\partial J^E}{\partial s}, X\right] e^{-A^2}\right]\right). \tag{129}$$

*Here $\left[(J^E)^{-1} \frac{\partial J^E}{\partial s}, X\right]$ is an ordinary matrix commutator and not a supercommutator.*

**Pf.** : We can write

$$\frac{\partial J^E}{\partial s} = [J^E, R] \tag{130}$$

with

$$R = \frac{1}{2} (J^E)^{-1} \frac{\partial J^E}{\partial s}. \tag{131}$$

Now

$$\begin{aligned} \frac{\partial A'}{\partial s} &= 0, \\ \frac{\partial A''}{\partial s} &= \left[\frac{\partial A''}{\partial s}, (J^E)^{-1} \frac{\partial J^E}{\partial s}\right] = [A'', 2R]. \end{aligned} \tag{132}$$

Thus

$$\frac{\partial A}{\partial s} = [A'', R] = [A, R] + [X, R]. \tag{133}$$

As an overall conjugation does not change the trace, we may equally well assume that

$$\begin{aligned}\frac{\partial J^E}{\partial s} &= 0, \\ \frac{\partial A}{\partial s} &= [X, R].\end{aligned} \quad (134)$$

That is, in effect we have a fixed $\mathbf{Z}_2$-grading on $E_{\mathbf{C}}$ and we are varying the superconnection $A$. It is known in this case [30] that

$$\frac{\partial}{\partial s}\operatorname{ch}(A_{\mathbf{C}}) = d\left(-(2i\pi)^{-1/2}\,\phi\operatorname{tr}_s\left[\frac{\partial A}{\partial s}e^{-A^2}\right]\right). \quad (135)$$

Thus

$$\frac{\partial}{\partial s}p(A', J^E(s)) = d\left(-(2i\pi)^{-1/2}\frac{1}{\sqrt{\epsilon}}\,\phi\operatorname{tr}\left[J^E\,[X,R]\,e^{-A^2}\right]\right). \quad (136)$$

The proposition follows. ■

**Definition 17** *Given $(\cdot,\cdot)$ and $A'$, let $J_1^E$ and $J_2^E$ be two choices of $J^E$. Let $J^E(s)$ be a smooth 1-parameter family of $J^E$'s such that $J^E(0) = J_2^E$ and $J^E(1) = J_1^E$. Define $\widetilde{p}(A', J_1^E, J_2^E) \in \Omega^{4*-\epsilon}(M)/\operatorname{im}(d)$ by*

$$\widetilde{p}(A', J_1^E, J_2^E) = \int_0^1 ds\,(2i\pi)^{-1/2}\frac{1}{2\sqrt{\epsilon}}\phi\operatorname{tr}\left[J^E\left[(J^E)^{-1}\frac{\partial J^E}{\partial s}, X\right]e^{-A^2}\right]. \quad (137)$$

**Corollary 4** *One has*

$$d\widetilde{p}(A', J_1^E, J_2^E) = p(A', J_1^E) - p(A', J_2^E). \quad (138)$$

*In particular, the de Rham cohomology class $p(A') \in \mathrm{H}^{4*+1-\epsilon}(M; \mathbf{R})$ of $p(A', J^E)$ is independent of $J^E$. Also, $\widetilde{p}(A', J_1^E, J_2^E)$ is independent of the choice of the 1-parameter family $J^E(s)$ in its definition.*

**Pf.** : This follows from Proposition 18. ■

**Corollary 5** *Fix $E$, $\nabla^E$ and $(\cdot,\cdot)$ as in Subsection 3.1. Then the de Rham cohomology class $p(\nabla^E) \in H^{4*+1-\epsilon}(M;\mathbf{R})$ of $p(\nabla^E, J^E)$ is independent of $J^E$. More precisely, let $s \in \mathbf{R}$ parametrize a smooth family of automorphisms $J^E(s)$ of $E$ as in Subsection 3.1. If $\epsilon = 1$ then*

$$\frac{\partial}{\partial s} p(\nabla^E, J^E(s)) = d\left(-\frac{1}{4\pi} \operatorname{tr}\left[\frac{\partial J^E}{\partial s} \omega \sin\left(\frac{\omega^2}{8\pi}\right)\right]\right). \tag{139}$$

*If $\epsilon = -1$ then*

$$\frac{\partial}{\partial s} p(\nabla^E, J^E(s)) = d\left(-\frac{1}{4\pi} \operatorname{tr}\left[\frac{\partial J^E}{\partial s} \omega \cos\left(\frac{\omega^2}{8\pi}\right)\right]\right). \tag{140}$$

**Pf. :** Applying Proposition 18,

$$\begin{aligned}
\frac{\partial}{\partial s} p(\nabla^E, J^E(s)) &= d\left((2i\pi)^{-1/2} \frac{1}{4\sqrt{\epsilon}} \phi \operatorname{tr}\left[J^E \left[(J^E)^{-1}\frac{\partial J^E}{\partial s}, \omega\right] e^{\frac{\omega^2}{4}}\right]\right) \\
&= d\left(\frac{1}{8i\pi} \frac{1}{\sqrt{\epsilon}} \operatorname{tr}\left[J^E \left[(J^E)^{-1}\frac{\partial J^E}{\partial s}, \omega\right] e^{\frac{\omega^2}{8i\pi}}\right]\right) \\
&= d\left(\frac{1}{8i\pi} \frac{1}{\sqrt{\epsilon}} \operatorname{tr}\left[J^E \left((J^E)^{-1}\frac{\partial J^E}{\partial s} \omega\right.\right.\right. \\
&\qquad\qquad \left.\left.\left. - \omega (J^E)^{-1}\frac{\partial J^E}{\partial s}\right) e^{\frac{\omega^2}{8i\pi}}\right]\right) \\
&= d\left(\frac{1}{8i\pi} \frac{1}{\sqrt{\epsilon}} \operatorname{tr}\left[\left(J^E \left((J^E)^{-1}\frac{\partial J^E}{\partial s} \omega\right)\right.\right.\right. \\
&\qquad\qquad \left.\left.\left. + \omega J^E \left((J^E)^{-1}\frac{\partial J^E}{\partial s}\right)\right) e^{\frac{\omega^2}{8i\pi}}\right]\right) \\
&= d\left(\frac{1}{4i\pi} \frac{1}{\sqrt{\epsilon}} \operatorname{tr}\left[J^E \left((J^E)^{-1}\frac{\partial J^E}{\partial s}\right) \omega\, e^{\frac{\omega^2}{8i\pi}}\right]\right). \tag{141}
\end{aligned}$$

The corollary follows. ∎

### 3.3  $\overline{L}$-Groups and Real-Valued Eta-Invariants

**Definition 18** *The generators of $L^0_\epsilon(M)$ are the triples $\mathcal{E} = (E, \nabla^E, (\cdot,\cdot)_E)$ where*

- $E$ is a real vector bundle on $M$.

- $\nabla^E$ is a flat connection on $E$.

- $(\cdot, \cdot)_E$ is a nondegenerate $\epsilon$-symmetric bilinear form on $E$ which is covariantly-constant with respect to $\nabla^E$.

**Definition 19** *The $L^0_\epsilon(M)$-relations are given by*
*1. Three $L^0_\epsilon(M)$-generators $\mathcal{E}^1$, $\mathcal{E}^2$ and $\mathcal{E}^3$, along with a short exact sequence of flat vector bundles*

$$0 \longrightarrow E^1 \xrightarrow{i} E^2 \xrightarrow{j} E^3 \longrightarrow 0 \tag{142}$$

*such that $i$ and $j^T$ are isometries.*
*2. An $L^0_\epsilon(M)$-generator of the form $\mathcal{V} \oplus \mathcal{V}^T = \left(V \oplus V^T, \nabla^V \oplus (\nabla^V)^T, (\cdot, \cdot)_{V \oplus V^T}\right)$, where*

- $V$ is a real vector bundle on $M$.

- $\nabla^V$ is a flat connection on $V$.

- $(\cdot, \cdot)_{V \oplus V^T}$ is given by

$$(v \oplus v', w \oplus w')_{V \oplus V^T} = w'(v) + \epsilon v'(w). \tag{143}$$

Note that $j^T$ gives a canonical splitting of the sequence (142). Thus

$$\mathcal{E}^2 \cong \mathcal{E}^1 \oplus \mathcal{E}^3. \tag{144}$$

**Definition 20** *The group $L^0_\epsilon(M)$ is the quotient of the free abelian group generated by the $L^0_\epsilon(M)$-generators, by the subgroup generated by the $L^0_\epsilon(M)$-relations $\mathcal{E}^2 - \mathcal{E}^1 - \mathcal{E}^3$ and $\mathcal{V} \oplus \mathcal{V}^T$.*

**Example 3 :** One has $L^0_1(\text{pt.}) = \mathbf{Z}$, where the right-hand-side is represented by the signature of the symmetric form $(\cdot, \cdot)_E$. Also, $L^0_{-1}(S^1) = \mathbf{Z}[SO(2)]$, the group ring of $SO(2)$. To see this, given an $L^0_{-1}(S^1)$-generator $\mathcal{E} = (E, \nabla^E, (\cdot, \cdot)_E)$, let $H\left(\nabla^E\right)$ be the holonomy of $\nabla^E$ around the circle with respect to a basepoint. Then $H\left(\nabla^E\right)$ lies in $Sp(2m, \mathbf{R})$, where $\dim(E) = 2m$. We can put $H\left(\nabla^E\right)$ into a normal form consisting of a direct sum of $2 \times 2$ hyperbolic and elliptic factors. Then the element of $\mathbf{Z}[SO(2)]$ corresponding to $\mathcal{E}$ is the sum of the elliptic factors.

**Remark 4 :** In the case of integral rather than real structures, similar groups $L^{2k}(\pi, \mathbf{Z}, v)$ were defined in [24, p. 147], where $(-1)^k = \epsilon$.

Given the sequence (142), choose automorphisms $\{J^{E^i}\}_{i=1}^3$ on $\{E^i\}_{i=1}^3$. Then using the splitting (144), we can write

$$d\widetilde{p}(\nabla^{E^2}, J^{E^2}, J^{E^1} \oplus J^{E^3}) = \sum_{i=1}^{3} (-1)^i \, p(\nabla^{E^i}, J^{E^i}). \tag{145}$$

**Proposition 19** *The assignment of $p(\nabla^E)$ to $E$ extends to a map $p : L_\epsilon^0(M) \to \mathrm{H}^{4*+1-\epsilon}(M; \mathbf{R})$.*

**Pf. :** We have to show that $p$ vanishes on the $L_\epsilon^0(M)$-relations.
1. For a relation coming from a short exact sequence (142), this follows from (145).
2. For a relation $\mathcal{V} \oplus \mathcal{V}^T$, choose a positive-definite inner product $h^V$ on $V$. Let

$$\widehat{h}^V : V \to V^T \tag{146}$$

be the isomorphism defined by

$$(\widehat{h}^V v)(w) = h^V(v, w). \tag{147}$$

Define an automorphism $J_0^{V \oplus V^T}$ by

$$J_0^{V \oplus V^T}(v \oplus v') = \epsilon \, (\widehat{h}^V)^{-1} v' \oplus \widehat{h}^V v. \tag{148}$$

Then $\omega^2$ is diagonal with respect to the splitting $V \oplus V^T$, whereas $J_0^{V \oplus V^T}$ is off-diagonal. It follows that

$$p(\nabla^{V \oplus V^T}, J_0^{V \oplus V^T}) = 0. \tag{149}$$

Thus $p(\nabla^{V \oplus V^T}) = 0$. ∎

**Definition 21** *An $\widehat{L}_\epsilon^0(M)$-generator is a triple $\mathcal{E} = (E, J^E, \rho)$ where*

- $E$ is a real vector bundle on $M$ with a flat connection $\nabla^E$ and a nondegenerate $\epsilon$-symmetric bilinear form $(\cdot, \cdot)_E$ which is covariantly constant with respect to $\nabla^E$.

- $J^E$ is an automorphism of the topological vector bundle $E$ such that $(J^E)^2 = \epsilon$, $(J^E e, J^E e')_E = (e, e')_E$ and $< e, e' >_E = (e, J^E e')_E$ is a positive-definite inner product on $E$.

- $\rho \in \Omega^{4*-\epsilon}(M)/\mathrm{im}(d)$.

**Definition 22** *The $\widehat{L}^0_\epsilon(M)$-relations are given by*

- *Three $\widehat{L}^0_\epsilon(M)$-generators $\mathcal{E}^1$, $\mathcal{E}^2$ and $\mathcal{E}^3$, along with a short exact sequence (142) such that*

$$\rho_2 = \rho_1 + \rho_3 + \widetilde{p}(\nabla^{E^2}, J^{E^2}, J^{E^1} \oplus J^{E^3}). \tag{150}$$

- *An $\widehat{L}^0_\epsilon(M)$-generator $\mathcal{V} \oplus \mathcal{V}^T = (V \oplus V^T, J_0^{V \oplus V^T}, 0)$, where $J_0^{V \oplus V^T}$ is defined as in (148).*

**Definition 23** *The group $\widehat{L}^0_\epsilon(M)$ is the quotient of the free abelian group generated by the $\widehat{L}^0_\epsilon(M)$-generators, by the subgroup generated by the $\widehat{L}^0_\epsilon(M)$-relations $\mathcal{E}^2 - \mathcal{E}^1 - \mathcal{E}^3$ and $\mathcal{V} \oplus \mathcal{V}^T$.*

**Proposition 20** *The assignment of $p(\nabla^E, J^E) - d\rho$ to $(E, J^E, \rho)$ extends to a map $p' : \widehat{L}^0_\epsilon(M) \to \Omega^{4*+1-\epsilon}(M)$.*

**Pf.** : This follows from (145) and (149). ∎

**Definition 24** *Let $\overline{L}^0_\epsilon(M)$ be the kernel of $p'$.*

There is a complex

$$\mathrm{H}^{4*-\epsilon}(M; \mathbf{R}) \xrightarrow{a} \overline{L}^0_\epsilon(M) \xrightarrow{b} L^0_\epsilon(M) \xrightarrow{p} \mathrm{H}^{4*+1-\epsilon}(M; \mathbf{R}), \tag{151}$$

where $a(\sigma) = [(0, 0, \sigma)]$ and $b\left(\sum_i n_i \left[(E^i, J^{E^i}, \rho^i)\right]\right) = \sum_i n_i [E^i]$.

**Proposition 21** *The complex (151) is exact.*

**Pf.** : Exactnesss at $L^0_\epsilon(M)$ : Given $\sum_i n_i[F^i] \in L^0_\epsilon(M)$ such that $p(\sum_i n_i[F^i]) = 0$, choose automorphisms $\{J^{F^i}\}$ on the $\{F^i\}$. As $\sum_i n_i p(\nabla^{F^i}, J^{F^i}) \in \Omega^{even}(M)$ represents $p(\sum_i n_i[F^i])$ in de Rham cohomology, there is a $\tau \in \Omega^{odd}(M)/\mathrm{im}(d)$ such that $\sum_i n_i\, p(\nabla^{F^i}, J^{F^i}) = d\tau$. Then $[(0,0,\tau)] + \sum_i n_i[(F^i, J^{F^i}, 0)]$ is an element of $\overline{L}^0_\epsilon(M)$ whose image under $b$ is $\sum_i n_i[F^i]$.

Exactness at $\overline{L}^0_\epsilon(M)$ : Given $\sum_i n_i[(F^i, J^{F^i}, \rho^i)] \in \overline{L}^0_\epsilon(M)$ such that $\sum_i n_i[F^i]$ vanishes in $L^0_\epsilon(M)$, there are flat real vector bundles $\{G^j, H^j, I^j\}$ on $M$ and short exact sequences

$$0 \longrightarrow G^j \longrightarrow H^j \longrightarrow I^j \longrightarrow 0, \tag{152}$$

along with flat real vector bundles $\{V^k \oplus V^{kT}\}$ on $M$ such that

$$\sum_i n_i\, F^i = \sum_j m_j \left(G^j - H^j + I^j\right) + \sum_k l_k \left(V^k \oplus V^{kT}\right) \tag{153}$$

for some integers $\{m_j\}$ and $\{l_k\}$. We also have

$$\sum_i n_i\, p(\nabla^{F^i}, J^{F^i}) = d\sum_i n_i\, \rho^i \tag{154}$$

Choose automorphisms $\{J^{G^j}, J^{H^j}, J^{I^j}\}$ and $\{J^{V^k\oplus V^{kT}}\}$ so that (with a slight abuse of terminology) coincident terms in (153) have the same automorphisms. Then

$$\begin{aligned}\sum_i n_i\, p(\nabla^{F^i}, J^{F^i}) &= \sum_j m_j \left(p(\nabla^{G^j}, J^{G^j}) - p(\nabla^{H^j}, J^{H^j}) - p(\nabla^{I^j}, J^{I^j})\right) \\ &+ \sum_k l_k\, p(\nabla^{V^k\oplus V^{kT}}, J^{V^k\oplus V^{kT}}).\end{aligned} \tag{155}$$

Choose positive-definite inner products $\{h^{V^k}\}$ on $\{V^k\}$ and define $J_0^{V^k\oplus V^{kT}}$ as in (148). We have

$$d\widetilde{p}(\nabla^{H^j}, J^{H^j}, J^{G^j}\oplus J^{I^j}) = -p(\nabla^{G^j}, J^{G^j}) + p(\nabla^{H^j}, J^{H^j}) - p(\nabla^{I^j}, J^{I^j}) \tag{156}$$

and

$$d\widetilde{p}(\nabla^{V^k\oplus V^{kT}}, J_0^{V^k\oplus V^{kT}}, J^{V^k\oplus V^{kT}}) = -p(\nabla^{V^k\oplus V^{kT}}, J^{V^k\oplus V^{kT}}). \tag{157}$$

From (154), (155), (156) and (157), we have that

$$\begin{aligned}\tau \equiv & \sum_i n_i \rho^i + \sum_j m_j \, \widetilde{p}(\nabla^{H^j}, J^{H^j}, J^{G^j} \oplus J^{I^j}) \\ & + \sum_k l_k \, \widetilde{p}(\nabla^{V^k \oplus V^{kT}}, J_0^{V^k \oplus V^{kT}}, J^{V^k \oplus V^{kT}})\end{aligned} \tag{158}$$

lies in $\mathrm{H}^{even}(B; \mathbf{R})$. Then in $\overline{L}_\epsilon^0(M)$, we have

$$\begin{aligned}\sum_i n_i[(F^i, J^{F^i}, \rho^i)] = & \, a(\tau) + \sum_j m_j \left([(G^j, J^{G^j}, 0)] \right. \\ & \left. - [(H^j, J^{H^j}, \widetilde{p}(\nabla^{H^j}, J^{H^j}, J^{G^j} \oplus J^{I^j}))] + [(I^j, h^{I^j}, 0)]\right) \\ & + \sum_k l_k \left([(V^k \oplus V^{kT}, J^{V^k \oplus V^{kT}}, 0)] - [(V^k \oplus V^{kT}, \right. \\ & \left. J_0^{V^k \oplus V^{kT}}, \widetilde{p}(\nabla^{V^k \oplus V^{kT}}, J_0^{V^k \oplus V^{kT}}, J^{V^k \oplus V^{kT}}))]\right) \\ = & \, a(\tau). \end{aligned} \tag{159}$$

The proposition follows. ∎

**Proposition 22** *Define $i_0 : M \to [0,1] \times M$ by $i_0(m) = (0, m)$ and $i_1 : M \to [0,1] \times M$ by $i_1(m) = (1, m)$. Then for any $\widetilde{l} \in \overline{L}_\epsilon^0([0,1] \times M)$, one has $i_0^* \widetilde{l} = i_1^* \widetilde{l}$.*

**Pf. :** The proof is similar to that of Proposition 4. We omit the details. ∎

**Corollary 6** *Let $Z$ be a smooth connected manifold. Given a smooth map $\phi : Z \to M$ and an element $l \in \overline{L}_\epsilon^0(M)$, the pullback $\phi^* l \in \overline{L}_\epsilon^0(Z)$ only depends on the (smooth) homotopy class of $\phi$.*

**Corollary 7** *If $f : M \to M'$ is a (smooth) homotopy equivalence then it induces an isomorphism $f^* : \overline{L}_\epsilon^0(M') \to \overline{L}_\epsilon^0(M)$.*

We now show how to detect elements of $\overline{L}_\epsilon^0(M)$ by real-valued eta-invariants. Let $Z$ be a connected closed oriented $n$-dimensional smooth manifold. Let $g^{TZ}$ be a Riemannian metric on $Z$ and let $\nabla^{TZ}$ be the Levi-Civita connection on $TZ$. Consider the real vector space of smooth differential forms $\Omega(Z)$. Let $*$ denote the Hodge duality operator on $\Omega(Z)$. Let $(\cdot,\cdot)_Z$ be the nondegenerate bilinear form on $\Omega(Z)$ such that for $e \in \Omega^i(Z)$ and $e' \in \Omega^{n-i}(Z)$,

$$(e, e')_Z = (-1)^{\frac{i(i-1)}{2}} \int_Z e \wedge e'. \tag{160}$$

Define $J^Z \in \text{Aut}(\Omega(Z))$ so that for $e \in \Omega^i(Z)$,

$$J^Z(e) = (-1)^{\frac{i(i-1)}{2}} * e. \tag{161}$$

Let $<\cdot,\cdot>_Z$ be the standard inner product on $\Omega(Z)$, namely

$$<e, e'>_Z = \int_Z e \wedge *e'. \tag{162}$$

Let $d^Z$ be exterior differentiation on $\Omega(Z)$ and let $\left(d^Z\right)^T$ be its formal adjoint.

**Lemma 5** *Put $\epsilon_n = (-1)^{\frac{n(n-1)}{2}}$. Then*
*a. $(\cdot,\cdot)_Z$ is $\epsilon_n$-symmetric.*
*b. $\left(J^Z\right)^2 = \epsilon_n$.*
*c. $(J^Z e, J^Z e')_Z = (e, e')_Z$.*
*d. $<e, e'>_Z = (e, J^Z e')_Z$.*
*e. $\left(d^Z\right)^T = -\left(J^Z\right)^{-1} d^Z J^Z$.*

**Pf.** : Parts a.-d. follow from a straightforward calculation. To prove the last part, if $e \in \Omega^i(Z)$ and $e' \in \Omega^{i+1}(Z)$ then

$$\begin{aligned}
<e, \left(d^Z\right)^T e'>_Z &= <d^Z e, e'>_Z = (d^Z e, J^Z e')_Z \\
&= (-1)^{\frac{i(i+1)}{2}} \int_Z d^Z e \wedge J^Z e' = -(-1)^{\frac{i(i-1)}{2}} \int_Z e \wedge d^Z J^Z e' \\
&= -(e, d^Z J^Z e')_Z = -<e, \left(J^Z\right)^{-1} d^Z J^Z e'>_Z. \tag{163}
\end{aligned}$$

The lemma follows. ∎

Let $\mathcal{E} = (E, J^E, \rho)$ be an $\widehat{L}_\epsilon^0(Z)$-generator. Consider the real vector space $\Omega(Z; E)$. Put

$$(\cdot,\cdot)_{Z,E} = (\cdot,\cdot)_{\Omega(Z)} \otimes_{C^\infty(Z)} (\cdot,\cdot)_E, \quad J^{Z,E} = J^Z \otimes_{C^\infty(Z)} J^E. \tag{164}$$

Then $(\cdot,\cdot)_{Z,E}$ is $\epsilon\epsilon_n$-symmetric and $\left(J^{Z,E}\right)^2 = \epsilon\epsilon_n$. Put

$$<e,e'>_{Z,E} = (e, J^{Z,E}e')_{Z,E} \qquad (165)$$

Let $d^{Z,E}$ be exterior differentiation on $\Omega(Z;E)$ using the flat connection $\nabla^E$ on $E$. Then the formal adjoint $\left(d^{Z,E}\right)^T$ satisfies

$$\left(d^{Z,E}\right)^T = -\left(J^{Z,E}\right)^{-1} d^{Z,E} J^{Z,E}. \qquad (166)$$

We assume for the rest of this subsection that $n$ is odd and $\epsilon\epsilon_n = -1$, i.e. $n \equiv -\epsilon \pmod{4}$.

**Definition 25** *The operator $D$ on $\Omega(Z;E)$ is given by*

$$D = J^{Z,E}d^{Z,E} + d^{Z,E}J^{Z,E}. \qquad (167)$$

**Proposition 23** *The operator $D$ is a real self-adjoint elliptic operator which commutes with $J^{Z,E}$ and satisfies*

$$D^2 = \left(d^{Z,E}\right)^T d^{Z,E} + d^{Z,E}\left(d^{Z,E}\right)^T. \qquad (168)$$

*Its kernel is isomorphic to $H^*(Z;E)$.*

**Pf. :** It is easy to see that $D$ commutes with $J^{Z,E}$. Its adjoint is given by

$$\begin{aligned} D^T &= \left(d^{Z,E}\right)^T \left(J^{Z,E}\right)^T + \left(J^{Z,E}\right)^T \left(d^{Z,E}\right)^T \\ &= \left(-J^{Z,E}d^{Z,E}(J^{Z,E})^{-1}\right)\left(-J^{Z,E}\right) \\ &\quad + \left(-J^{Z,E}\right)\left(-J^{Z,E}d^{Z,E}(J^{Z,E})^{-1}\right) \\ &= J^{Z,E}d^{Z,E} + d^{Z,E}J^{Z,E}. \end{aligned} \qquad (169)$$

Finally,

$$\begin{aligned} D^2 &= J^{Z,E}d^{Z,E}J^{Z,E}d^{Z,E} + d^{Z,E}J^{Z,E}d^{Z,E}J^{Z,E} \\ &= \left(d^{Z,E}\right)^T d^{Z,E} + d^{Z,E}\left(d^{Z,E}\right)^T. \end{aligned} \qquad (170)$$

Thus $D$ is elliptic. As $\operatorname{Ker}(D^2) \cong H^*(Z;E)$, the proposition follows. ∎

In particular, $\dim(\operatorname{Ker}(D))$ is independent of $g^{TZ}$ and $J^E$.

Let $\eta(D) \in \mathbf{R}$ be the eta-invariant of Atiyah-Patodi-Singer [2]. Let $L(TZ, \nabla^{TZ}) \in \Omega^{4*}(Z)$ be the Hirzebruch signature integrand.

**Definition 26** *Put $\eta(\mathcal{E}) = \frac{1}{2}\eta(D) - \int_Z L(TZ, \nabla^{TZ}) \wedge \rho$.*

**Proposition 24** *There is an extension of $\eta(\mathcal{E})$ to a linear map $\eta : \widehat{L}^0_\epsilon(Z) \to \mathbf{R}$.*

**Pf. :** We must show that $\eta$ vanishes on the $\widehat{L}^0_\epsilon(Z)$-relations. First, given an $\widehat{L}^0_\epsilon(Z)$-relation as in (150), we must show that

$$\frac{1}{2}\eta(D_2) - \frac{1}{2}\eta(D_1) - \frac{1}{2}\eta(D_3) = \int_Z L(TZ, \nabla^{TZ}) \wedge \widetilde{p}\left(\nabla^{E^2}, J^{E^2}, J^{E^1} \oplus J^{E^3}\right). \tag{171}$$

Equivalently, if $(E, \nabla^E, (\cdot, \cdot))$ is a generator of $L^0_\epsilon(Z)$ and $J^E_1$, $J^E_2$ are two choices of $J^E$, let $D_1$ and $D_2$ be the corresponding self-adjoint operators. Then we must show that

$$\frac{1}{2}\eta(D_1) - \frac{1}{2}\eta(D_2) = \int_Z L(TZ, \nabla^{TZ}) \wedge \widetilde{p}\left(\nabla^E, J^E_1, J^E_2\right). \tag{172}$$

Let $s \in [0,1]$ parametrize a smooth family of automorphisms $J^E(s)$ of $E$ such that $J^E(0) = J^E_2$, $J^E(1) = J^E_1$ and $J^E(s)$ is constant near the endpoints of $[0,1]$. Put $\widetilde{Z} = [0,1] \times Z$. Let $g^{T\widetilde{Z}}$ be the product metric on $\widetilde{Z}$. Let $\pi : \widetilde{Z} \to Z$ be projection onto the second factor. Put $\widetilde{E} = \pi^* E$ and $\nabla^{\widetilde{E}} = \pi^* \nabla^E$. We obtain an automorphism $J^{\widetilde{E}}$ of $\widetilde{E}$ by piecing together the $J^E(s)$'s. Define the operator

$$\widetilde{D} = d^{\widetilde{Z}, \widetilde{E}} + \left(d^{\widetilde{Z}, \widetilde{E}}\right)^T \tag{173}$$

on $\widetilde{Z}$, with the Atiyah-Patodi-Singer boundary conditions [2]. Then $\widetilde{D}$ anti-commutes with $J^{\widetilde{E}}$ and we can talk about its index. From the Atiyah-Patodi-Singer index theorem [2],

$$\mathrm{ind}(\widetilde{D}) = \int_{\widetilde{Z}} \text{char. form }(g^{T\widetilde{Z}}, J^{\widetilde{E}}) - \left(\frac{1}{2}\eta(D_1) - \frac{1}{2}\eta(D_2)\right), \tag{174}$$

where char. form $(g^{T\widetilde{Z}}, J^{\widetilde{E}})$ is the characteristic form of the operator $\widetilde{D}$. To compute this form, let $\{e_i\}_{i=1}^{n+1}$ be a local oriented orthonormal basis of $T\widetilde{Z}$. Using the Riemannian metric, we can identify $T\widetilde{Z}$ and its dual. If $V$ is a vector on $\widetilde{Z}$, put

$$\begin{aligned} c(V) &= (V\wedge) - i_V, \\ \widehat{c}(V) &= (V\wedge) + i_V. \end{aligned} \tag{175}$$

Then locally,

$$
\begin{aligned}
\widetilde{D} &= \sum_{j=1}^{n+1} \left( (e_j \wedge) \nabla^{\widetilde{E}}_{e_j} - i_{e_j} (\nabla^{\widetilde{E}})^T_{e_j} \right) \\
&= \sum_{j=1}^{n+1} \frac{1}{2} \left( (e_j \wedge) - i_{e_j} \right) \left( \nabla^{\widetilde{E}}_{e_j} + (\nabla^{\widetilde{E}})^T_{e_j} \right) \\
&\quad + \sum_{j=1}^{n+1} \frac{1}{2} \left( (e_j \wedge) + i_{e_j} \right) \left( \nabla^{\widetilde{E}}_{e_j} - (\nabla^{\widetilde{E}})^T_{e_j} \right) \\
&= \sum_{j=1}^{n+1} \left( c(e_j) \nabla^{\widetilde{E},u}_{e_j} - \frac{1}{2} \widehat{c}(e_j)\, \omega_j(\nabla^{\widetilde{E}}, J^{\widetilde{E}}) \right). \quad (176)
\end{aligned}
$$

We can now use the Getzler rescaling to compute the index density [5], in which

$$
\begin{aligned}
\partial_j &\to u^{-\frac{1}{2}} \partial_j, \\
c_j &\to u^{-\frac{1}{2}} (e_j \wedge) - u^{\frac{1}{2}} i_{e_j}, \\
\widehat{c}_j &\to \widehat{c}_j.
\end{aligned}
\quad (177)
$$

The result is

$$
\text{char. form}\,(g^{T\widetilde{Z}}, J^{\widetilde{E}}) = L(T\widetilde{Z}, \nabla^{T\widetilde{Z}}) \wedge p(\nabla^{\widetilde{E}}, J^{\widetilde{E}}). \quad (178)
$$

In order to compute $\mathrm{ind}(\widetilde{D})$, one can homotop to $J^E(s)$ being constant in $s$, in which case one can check that $\mathrm{ind}(\widetilde{D}) = 0$. Then (174) becomes

$$
\begin{aligned}
\frac{1}{2}\eta(D_1) - \frac{1}{2}\eta(D_2) &= \int_{\widetilde{Z}} L(T\widetilde{Z}, \nabla^{T\widetilde{Z}}) \wedge p(\nabla^{\widetilde{E}}, J^{\widetilde{E}}) \\
&= \int_Z L(TZ, \nabla^{TZ}) \wedge \widetilde{p}\left(\nabla^E, J_1^E, J_2^E\right). \quad (179)
\end{aligned}
$$

Next, supppose that we have a $\widehat{L}^0_\epsilon(Z)$-relation $\mathcal{V} \oplus \mathcal{V}^T = (V \oplus V^T, J_0^{V \oplus V^T}, 0)$. Then with respect to the splitting $V \oplus V^T$, the operator $D$ anticommutes with the matrix $\begin{pmatrix} I & 0 \\ 0 & -I \end{pmatrix}$. It follows that the spectrum of $D$ is symmetric around the origin and so $\eta(D) = 0$.

The proposition follows. ∎

**Proposition 25** *The restriction of $\eta : \widehat{L}^0_\epsilon(Z) \to \mathbf{R}$ to $\overline{L}^0_\epsilon(Z)$ is independent of $g^{TZ}$.*

**Pf.** : Let $g_1^{TZ}$ and $g_2^{TZ}$ be two Riemannian metrics on $Z$. Let $\{\mathcal{E}^j\} = \{(E^j, J^{E^j}, \rho^j)\}$ be a finite set of $\widehat{L}^0_\epsilon(Z)$-generators. If $\sum_j n_j \, \mathcal{E}^j$ lies in $\overline{L}^0_\epsilon(Z)$ then
$$\sum_j n_j \, p(\nabla^{E^j}, J^{E^j}) = \sum_j n_j \, d\rho_j. \tag{180}$$

Let $D_1^j$ and $D_2^j$ be the corresponding self-adjoint operators. Then we must show that
$$\sum_j n_j \left( \frac{1}{2} \eta(D_1^j) - \frac{1}{2} \eta(D_2^j) \right) = \int_Z \left( L(TZ, \nabla_1^{TZ}) - L(TZ, \nabla_2^{TZ}) \right) \wedge \sum_j n_j \, \rho^j. \tag{181}$$

Let $\widetilde{Z}$, $\pi$, $\widetilde{E^j}$ and $\widetilde{\nabla^{E^j}}$ be as in the proof of Proposition 24. Put $\widetilde{J^{E^j}} = \pi^* J^{E^j}$. Let $s \in [0,1]$ parametrize a smooth family of Riemannian metrics $g^{TZ}(s)$ on $Z$ such that $g^{TZ}(0) = g_2^{TZ}$, $g^{TZ}(1) = g_1^{TZ}$ and $g^{TZ}(s)$ is constant near the endpoints of $[0,1]$. Put $g^{T\widetilde{Z}} = ds^2 + g^{TZ}(s)$, which is a product metric near the boundary. As in the proof of Proposition 24, we obtain

$$\begin{aligned}
\sum_j n_j \left( \frac{1}{2} \eta(D_1^j) - \frac{1}{2} \eta(D_2^j) \right) &= \int_{\widetilde{Z}} L(T\widetilde{Z}, \nabla^{T\widetilde{Z}}) \wedge \sum_j n_j \, p(\widetilde{\nabla^{E^j}}, \widetilde{J^{E^j}}) \\
&= \int_{\widetilde{Z}} L(T\widetilde{Z}, \nabla^{T\widetilde{Z}}) \wedge \pi^* \sum_j n_j \, d\rho^j \\
&= \int_{\widetilde{Z}} d \left( L(T\widetilde{Z}, \nabla^{T\widetilde{Z}}) \wedge \sum_j n_j \, \pi^* \rho^j \right) \tag{182} \\
&= \int_Z \left( L(TZ, \nabla_1^{TZ}) - L(TZ, \nabla_2^{TZ}) \right) \wedge \sum_j n_j \, \rho^j.
\end{aligned}$$

The proposition follows. ∎

**Corollary 8** *Given $\epsilon = \pm 1$, let $Z$ be a connected closed oriented $n$-dimensional smooth manifold with $n \equiv -\epsilon \pmod{4}$. Let $M$ be a smooth connected manifold and let $[Z, M]$ denote the (smooth) homotopy classes of maps from $Z$ to $M$. Then there is a well-defined pairing $[Z, M] \times \overline{L}^0_\epsilon(M) \to \mathbf{R}$ given by $(\phi, l) \to \eta(\phi^* l)$.*

**Pf. :** This follows from Corollary 6 and Proposition 25. ∎

For an important special case, let $E$ be a flat real vector bundle on $M$ with holonomy in $O(m) \times O(m)$ (or $U(m)$). Let $J^E$ be the covariantly-constant automorphism of $E$ given by the matrix $\begin{pmatrix} I_m & 0 \\ 0 & -I_m \end{pmatrix}$ (or the complex structure). Define $l \in \overline{L}_\epsilon^0(M)$ by $l = [(E, J^E, 0)]$. Then $\eta(\phi^* l)$ is the rho-invariant of the tangential signature operator, as considered in [2]. In this case the (mod $\mathbf{Z}$) reduction of the pairing $(\phi, l) \to \eta(\phi^* l)$ has a simple topological description. Namely, the signature operator of $Z$ determines a class $[\sigma_Z] \in K_{-1}(Z)$. We can push it forward to $\phi_*([\sigma_Z]) \in K_{-1}(M)$. There is a forgetful map $f$ taking $l$ to $f(l) \in K^{-1}(M; \mathbf{R}/\mathbf{Z})$, where the latter is as defined in [22]. Then the (mod $\mathbf{Z}$) reduction of $\eta(\phi^* l)$ comes from the pairing

$$K_{-1}(M) \times K^{-1}(M; \mathbf{R}/\mathbf{Z}) \to \mathbf{R}/\mathbf{Z}, \tag{183}$$

applied to $\phi_*([\sigma_Z])$ and $f(l)$ [22]. In particular, the (mod $\mathbf{Z}$) reduction is a cobordism invariant of $\phi$. One can see in examples that the unreduced pairing $(\phi, l) \to \eta(\phi^* l)$ is not a cobordism invariant of $\phi$. We do not know of a simple topological description of the unreduced pairing.

### 3.4 Number Operators

We use the notation of Subsection 3.2.

**Definition 27** *The triple $(E, (\cdot, \cdot)_E, J^E)$ is $\mathbf{Z}$-graded if*

- *The vector bundle $E$ is $\mathbf{Z}$-graded as $E = \bigoplus_{i=0}^n E^i$, with the number operator $N \in \mathrm{End}(E)$ acting on $E^i$ as multiplication by $i$.*

- *If $e \in E^i$ and $e' \in E^{i'}$ with $i + i' \neq n$ then $(e, e')_E = 0$.*

- *$J^E E^i = E^{n-i}$.*

*A flat pair $(A, X)$ is of degree $1$ if we can write*

$$A - X = A' = \sum_{j \in \mathbf{N}} A'_j, \tag{184}$$

*where $A'_1$ is a connection on $E$ which preserves the $\mathbf{Z}$-grading and $A'_j \in \Omega^j(M; \mathrm{Hom}(E^\bullet, E^{\bullet+1-j}))$ for $j \neq 1$.*

**Warning :** The $\mathbf{Z}_2$-grading on $E_{\mathbf{C}}$ given by $\frac{1}{\sqrt{\epsilon}} J^E$ does not come from the $\mathbf{Z}$-grading on $E$ given by $N$. The only compatibility between the two is that

$$J^E(N - \frac{n}{2}) + (N - \frac{n}{2})J^E = 0, \tag{185}$$

showing that $N - \frac{n}{2}$ is an odd operator.

In the rest of this subsection, we assume that $(E, (\cdot, \cdot)_E, J^E)$ is $\mathbf{Z}$-graded and that $(A, X)$ is a flat pair of degree 1. Recall the definition of $A'$ and $A''$ from (112).

**Proposition 26** *With respect to the $\mathbf{Z}$-grading on $E$ coming from $N$, $A'$ and $A''$ are flat superconnections, with $A'$ being of total degree 1 in the sense of Definition 41.*

**Pf. :** The flatness of $(A, X)$ is equivalent to

$$\begin{aligned} 0 &= A^2 + (X\sigma)^2, \\ 0 &= AX\sigma + X\sigma A. \end{aligned} \tag{186}$$

As in the proof of Proposition 16, write

$$A = \sum_{j \geq 0} A_j, \quad X = \sum_{j \geq 0} X_j \tag{187}$$

and

$$A_j = \sum_r \omega_{j,r} A_{j,r}, \quad X_j = \sum_r \omega_{j,r} X_{j,r}. \tag{188}$$

With respect to the $\mathbf{Z}_2$-grading coming from $J^E$, we obtain

$$\begin{aligned} 0 &= A^2 + (X\sigma)^2 \\ &= \sum_{j,k,r,s} \left( \omega_{j,r} A_{j,r} \omega_{k,s} A_{k,s} + \omega_{j,r} X_{j,r} \sigma \omega_{k,s} X_{k,s} \sigma \right) \\ &= \sum_{j,k,r,s} \left( \omega_{j,r} A_{j,r} \omega_{k,s} A_{k,s} + (-1)^k \omega_{j,r} X_{j,r} \omega_{k,s} X_{k,s} \right) \\ &= \sum_{j,k,r,s} (-1)^{k(j-1)} \omega_{j,r} \wedge \omega_{k,s} \left( A_{j,r} A_{k,s} + X_{j,r} X_{k,s} \right). \end{aligned} \tag{189}$$

Similarly,

$$\begin{aligned}
0 &= AX\sigma + X\sigma A \\
&= \sum_{j,k,r,s} \left(\omega_{j,r}\, A_{j,r}\, \omega_{k,s}\, X_{k,s}\, \sigma + \omega_{j,r}\, X_{j,r}\, \sigma\, \omega_{k,s}\, A_{k,s}\right) \\
&= \sum_{j,k,r,s} \left(\omega_{j,r}\, A_{j,r}\, \omega_{k,s}\, X_{k,s} + (-1)^k\, \omega_{j,r}\, X_{j,r}\, \omega_{k,s}\, A_{k,s}\right)\sigma \\
&= \sum_{j,k,r,s} (-1)^{k(j-1)}\, \omega_{j,r} \wedge \omega_{k,s}\, (A_{j,r}\, X_{k,s} + X_{j,r}\, A_{k,s})\,\sigma
\end{aligned}$$
(190)

Thus
$$0 = \sum_{j,k,r,s} (-1)^{k(j-1)}\, \omega_{j,r} \wedge \omega_{k,s}(A_{j,r} \pm X_{j,r})\,(A_{k,s} \pm X_{k,s}). \tag{191}$$

On the other hand, with respect to the $\mathbf{Z}_2$-grading coming from $N$, this implies

$$\begin{aligned}
0 &= \sum_{j,k,r,s} \omega_{j,r}\,(A_{j,r} \pm X_{j,r})\, \omega_{k,s}\,(A_{k,s} \pm X_{k,s}) \\
&= (A \pm X)^2.
\end{aligned}$$
(192)

Thus $A'$ and $A''$ are flat superconnections. It is clear that $A'$ is of total degree 1. ∎

Let us write $v = A'_0$ and $\nabla^E = A'_1$. The flatness of $A'$ implies that $v^2 = \nabla^E v = 0$. Thus we have a cochain complex of vector bundles

$$(E, v) : 0 \to E^0 \xrightarrow{v} E^1 \xrightarrow{v} \cdots \xrightarrow{v} E^n \to 0 \tag{193}$$

whose differential $v$ is covariantly constant with respect to $\nabla^E$.

**Definition 28** *For $m \in M$, let $H(E,v)_m = \bigoplus_{i=0}^n \mathrm{H}^i(E,v)_m$ be the cohomology of the complex $(E,v)_m$.*

As in [9, Section 2a], there is a $\mathbf{Z}$-graded vector bundle $H(E,v)$ on $M$ whose fiber over $m \in M$ is $H(E,v)_m$. Furthermore, there is a natural flat connection $\nabla^H$ on $H(E,v)$ which can be described as follows. Let $\psi : \mathrm{Ker}(v) \to H(E,v)$ be the quotient map. Let $s$ be a smooth section of $H(E,v)$. Let $e$ be a smooth section of $\mathrm{Ker}(v)$ such that $\psi(e) = s$. Then if $U$ is a vector field on $M$,
$$\nabla^H_U s = \psi\left(\nabla^E_U e\right). \tag{194}$$

**Definition 29** *Given smooth sections $s$ and $s'$ of $H(E, v)$, choose smooth sections $e$ and $e'$ of $\mathrm{Ker}(v)$ such that $\psi(e) = s$ and $\psi(e') = s'$. Define a bilinear form $(\cdot, \cdot)_H$ on $H(E, v)$ by*

$$(s, s')_H = (e, e')_E. \tag{195}$$

We now assume that $(A, X)$ is $(\cdot, \cdot)_E$-compatible.

**Proposition 27** *The form $(\cdot, \cdot)_H$ is well-defined. It is covariantly-constant with respect to $\nabla^H$.*

**Pf.** : It follows from (111) that if $e, e' \in C^\infty(M; E)$ then

$$0 = (ve, e') + (e, ve'). \tag{196}$$

Suppose that $e_1$ and $e_2$ are smooth sections of $\mathrm{Ker}(v)$ such that $\psi(e_1) = \psi(e_2)$. Then $e_1 - e_2$ can be written as $v(f)$ for some smooth section $f$ of $H(E, v)$. Thus

$$(e_1, e')_E - (e_2, e')_E = (v(f), e')_E = -(f, v(e'))_E = 0, \tag{197}$$

showing that $(\cdot, \cdot)_H$ is well-defined. Furthermore, it follows from (111) that if $U$ is a vector field on $M$ and $e, e' \in C^\infty(M; E)$ then

$$U(e, e')_E = (\nabla^E_U e, e')_E + (e, \nabla^E_U e')_E. \tag{198}$$

If $s$ and $s'$ are smooth sections of $H(E, v)$ and $e$, $e'$ are as above then

$$U(s, s')_H = U(e, e')_E = (\nabla^E_U e, e')_E + (e, \nabla^E_U e')_E = (\nabla^H_U s, s')_H + (s, \nabla^H_U s')_H. \tag{199}$$

The proposition follows. ∎

We have that $A_0^2 = \frac{1}{4}(vv^* + v^*v)$. By Hodge theory, $H(E, v) \cong \mathrm{Ker}(A_0^2)$. Also, $J^E$ commutes with $A_0^2$.

**Definition 30** *Let $J^H$ be the restriction of $J^E$ to $\mathrm{Ker}(A_0^2)$. Define an inner product on $H(E, v)$ by $< s, s' >_H = (s, J^H s')_H$.*

It is clear that $(J^H)^2 = \epsilon$, $(s, s')_H = (J^H s, J^H s')_H$ and that $< \cdot, \cdot >_H$ is positive-definite. In particular, it follows that $(\cdot, \cdot)_H$ is nondegenerate and so the triple $(H(E, v), \nabla^H, (\cdot, \cdot)_H)$ defines a flat duality bundle.

We now introduce a rescaling of $J^E$. Let us denote $J^E$ by $J^E(1)$ and for $t \in \mathbf{R}^+$, put $J^E(t) = J^E(1)\, t^{N-\frac{n}{2}}$. Then
$$\left(J^E(t)\right)^2 = J^E(1)\, t^{N-\frac{n}{2}}\, J^E(1)\, t^{N-\frac{n}{2}} = J^E(1)\, J^E(1)\, t^{-N+\frac{n}{2}}\, t^{N-\frac{n}{2}} = \epsilon. \quad (200)$$

Fixing $A'$, let $A_t$ and $X_t$ be the odd and even parts of $A'$ with respect to $J^E(t)$. Put
$$\begin{aligned} C_t &= t^{N/2}\, A_t\, t^{-N/2}, \\ D_t &= t^{N/2}\, X_t\, t^{-N/2}. \end{aligned} \quad (201)$$

Then one can check that the superconnection $C_t$ is symmetric with respect to $J^E(1)$. Explicitly,
$$\begin{aligned} C_t &= \sum_{j\geq 0} t^{\frac{1-j}{2}}\, A_{1,j}, \\ D_t &= \sum_{j\geq 0} t^{\frac{1-j}{2}}\, X_{1,j}. \end{aligned} \quad (202)$$

One has
$$\begin{aligned} p(A', J^E(t)) &= \frac{1}{\sqrt{\epsilon}}\, \phi\, \mathrm{tr}\left[J^E(t)\, e^{-A_t^2}\right] \\ &= \frac{1}{\sqrt{\epsilon}}\, \phi\, \mathrm{tr}\left[J^E(1)\, t^{N-\frac{n}{2}}\, e^{-A_t^2}\right] \\ &= \frac{1}{\sqrt{\epsilon}}\, \phi\, \mathrm{tr}\left[t^{-\frac{N}{2}+\frac{n}{4}}\, J^E(1)\, t^{\frac{N}{2}-\frac{n}{4}}\, e^{-A_t^2}\right] \\ &= \frac{1}{\sqrt{\epsilon}}\, \phi\, \mathrm{tr}\left[J^E(1)\, t^{\frac{N}{2}-\frac{n}{4}}\, e^{-A_t^2}\, t^{-\frac{N}{2}+\frac{n}{4}}\right] \\ &= \frac{1}{\sqrt{\epsilon}}\, \phi\, \mathrm{tr}\left[J^E(1)\, e^{-C_t^2}\right]. \end{aligned} \quad (203)$$

From Proposition 18,
$$\begin{aligned} \frac{\partial}{\partial t} p(A', J^E(t)) &= d\left((2i\pi)^{-1/2}\, \frac{1}{2t\sqrt{\epsilon}}\, \phi\, \mathrm{tr}\left[J^E(t)\left[N - \frac{n}{2}, X_t\right] e^{-A_t^2}\right]\right) \\ &= d\left((2i\pi)^{-1/2}\, \frac{1}{2t\sqrt{\epsilon}}\, \phi\, \mathrm{tr}\left[J^E(t)\, [N, X_t]\, e^{-A_t^2}\right]\right) \\ &= d\left((2i\pi)^{-1/2}\, \frac{1}{2t\sqrt{\epsilon}}\, \phi\, \mathrm{tr}\left[J^E(1)\, [N, D_t]\, e^{-C_t^2}\right]\right) \quad (204) \end{aligned}$$

**Definition 31** *Define $\widetilde{\eta}(t) \in \Omega^{4*-\epsilon}(M)$ by*

$$\widetilde{\eta}(t) = (2i\pi)^{-1/2} \frac{1}{2t\sqrt{\epsilon}} \phi \operatorname{tr}\left[ J^E(1) [N, D_t] e^{-C_t^2} \right] \tag{205}$$

From (204),

$$\frac{\partial}{\partial t} p(A', J^E(t)) = d\widetilde{\eta}(t). \tag{206}$$

**Proposition 28** *As $t \to \infty$,*

$$\begin{aligned} p(A', J^E(t)) &= p(\nabla^H, J^H) + O(t^{-1/2}), \\ \widetilde{\eta}(t) &= O(t^{-3/2}). \end{aligned} \tag{207}$$

**Pf. :** The proof is similar to that of [9, Theorem 2.13]. We omit the details. ■

**Corollary 9** *As elements of $\mathrm{H}^{4*+1-\epsilon}(M; \mathbf{R})$,*

$$p(A') = p(\nabla^H). \tag{208}$$

**Pf. :** This follows from (206) and Proposition 28. ■

Now consider the special case in which $A'_j = 0$ for $j > 1$. That is, the differential $v$ is covariantly-constant with respect to a flat connection $\nabla^E$.

**Proposition 29** *As $t \to 0$,*

$$\begin{aligned} p(A', J^E(t)) &= p(\nabla^E, J^E) + O(t), \\ \widetilde{\eta}(t) &= O(1). \end{aligned} \tag{209}$$

**Pf. :** In this case, we have

$$\begin{aligned} C_t &= \frac{\sqrt{t}}{2}(v^* + v) + \nabla^{E,u}, \\ D_t &= \frac{\sqrt{t}}{2}(v^* - v) + \frac{\omega}{2}. \end{aligned} \tag{210}$$

Using the fact that $p(A', J^E(t))$ is an even form and expanding in $t$, we find

$$\begin{aligned}
p(A', J^E(t)) &= \frac{1}{\sqrt{\epsilon}} \phi \operatorname{tr}\left[ J^E(1)\, e^{-C_t^2} \right] \\
&= \frac{1}{\sqrt{\epsilon}} \phi \operatorname{tr}\left[ J^E(1)\, e^{-(\nabla^{E,u})^2} \right] + O(t) \\
&= p(\nabla^E, J^E) + O(t).
\end{aligned} \quad (211)$$

Similarly,

$$\begin{aligned}
\widetilde{\eta}(t) &= (2i\pi)^{-1/2} \frac{1}{2t\sqrt{\epsilon}} \phi \operatorname{tr}\left[ J^E(1)\, [N, D_t]\, e^{-C_t^2} \right] \\
&= (2i\pi)^{-1/2} \frac{1}{4\sqrt{t\epsilon}} \phi \operatorname{tr}\left[ J^E(1)\, [N, v^* - v]\, e^{-C_t^2} \right] \\
&= O(1).
\end{aligned} \quad (212)$$

The proposition follows. ∎

**Definition 32** *Define $\widetilde{\eta} \in \Omega^{4*-\epsilon}(M)$ by*

$$\widetilde{\eta} = -\int_0^\infty \widetilde{\eta}(t)\, dt. \quad (213)$$

By Propositions 28 and 29, $\widetilde{\eta}$ is well-defined and satisfies

$$d\widetilde{\eta} = p\left(\nabla^E, J^E\right) - p\left(\nabla^H, J^H\right). \quad (214)$$

The eta-form $\widetilde{\eta}$ is a special case of that defined in [7, Section 2].

## 3.5 Fiber Bundles

Let $Z \to M \xrightarrow{\pi} B$ be a smooth fiber bundle with connected base $B$ and connected closed fibers $Z_b = \pi^{-1}(b)$ of dimension $n$. Let $TZ$ be the vertical tangent bundle of the fiber bundle and let $T^*Z$ be its dual bundle. We assume that $TZ$ is oriented. Let $L(TZ) \in \mathrm{H}^{4*}(M)$ be the Hirzebruch $L$-class of $TZ$.

Equip the fiber bundle with a horizontal distribution $T^H M$. Let $\Omega(Z)$ denote the infinite-dimensional real vector bundle on $B$ whose fiber over $b \in B$ is isomorphic to $\Omega(Z_b)$. Then

$$C^\infty(B; \Omega(Z)) \simeq C^\infty(M; \Lambda(T^*Z)) \quad (215)$$

and there is an isomorphism of real vector spaces

$$\Omega(M) \simeq \Omega(B; \Omega(Z)). \tag{216}$$

Define a nondegenerate $\epsilon_n$-symmetric bilinear form $(\cdot, \cdot)_Z$ on $\Omega(Z)$ as in (160). Let $N$ be the number operator of $\Omega(Z)$; it acts as multiplication by $j$ on $C^\infty(M; \Lambda^j(T^*Z))$.

Let $g^{TZ}$ be a vertical Riemannian metric on the fiber bundle. For notation, we let lower case Greek indices refer to horizontal directions, lower case Roman indices refer to vertical directions and upper case Roman indices refer to either. We let $\{\tau^J\}$ denote a local basis of 1-forms on $M$, with dual basis $\{e_J\}$ of tangent vectors. We will always take $\{e_i\}_{i=1}^n$ to be an oriented orthonormal framing of $TZ$. We will assume that the forms $\{\tau^\alpha\}$ are pulled back from a local basis of 1-forms on $B$, which we will also denote by $\{\tau^\alpha\}$. Exterior multiplication by a form $\phi$ will be denoted by $\phi\wedge$ and interior multiplication by a vector $v$ will be denoted by $i_v$. Using the horizontal distribution and vertical Riemannian metric, we can identify vertical vectors and vertical 1-forms. Exterior multiplication by $\tau^J$ will be denoted by $E^J$ and interior multiplication by $e_J$ will be denoted by $I^J$. We have that $E^j I^k + I^k E^j = \delta^{jk}$. If $X$ is a vertical vector (or 1-form), put

$$\begin{aligned} c(X) &= (X\wedge) - i_X, \\ \widehat{c}(X) &= (X\wedge) + i_X. \end{aligned} \tag{217}$$

Put $c^i = c(e_i)$ and $\widehat{c}^i = \widehat{c}(e_i)$.

In calculations we will sometimes assume that $B$ has a Riemannian metric $g^{TB}$ and $M$ has the Riemannian metric $g^{TM} = g^{TZ} \oplus \pi^* g^{TB}$, although all final results will be independent of $g^{TB}$. Let $\nabla^{TM}$ denote the corresponding Levi-Civita connection on $M$ and put $\nabla^{TZ} = P^{TZ}\nabla^{TM}$, a connection on $TZ$. As shown in [6, Theorem 1.9], $\nabla^{TZ}$ is independent of the choice of $g^{TB}$. The restriction of $\nabla^{TZ}$ to a fiber coincides with the Levi-Civita connection of the fiber. We will also denote by $\nabla^{TZ}$ the extension to a connection on $\Lambda(T^*Z)$.

We will use the Einstein summation convention freely, and write

$$\omega_{IJK} = \tau^I(\nabla^{TM}_{e_K} e_J). \tag{218}$$

As there is a vertical metric, we may raise and lower vertical indices freely.

The fundamental geometric tensors of the fiber bundle are its curvature, a $TZ$-valued horizontal 2-form on $M$, and the second fundamental form of the fibers, a $(T^H M)^*$-valued vertical symmetric form on $M$. The curvature 2-form $T$ is given in terms of the local framing by

$$T_j = -\omega_{\alpha\beta j} \tau^\alpha \wedge \tau^\beta. \tag{219}$$

Define a horizontal 1-form $k$ on $M$ by

$$k = \omega_{j\alpha j} \tau^\alpha, \tag{220}$$

the mean curvature 1-form to the fibers.

Define $J^Z \in \mathrm{Aut}(\Omega(Z))$ as in (161). Then $\left(\Omega(Z), (\cdot,\cdot)_Z, J^Z\right)$ is **Z**-graded by the number operator $N$. Let $d^M$ be exterior differentiation on $\Omega(M)$. Let $\left(d^M\right)^T$ be its adjoint, as considered in [9, Proposition 3.7]. Put

$$\begin{aligned}
A_Z &= \frac{1}{2}\left(\left(d^M\right)^T + d^M\right), \\
X_Z &= \frac{1}{2}\left(\left(d^M\right)^T - d^M\right).
\end{aligned} \tag{221}$$

**Proposition 30** *One has that $A_Z$ is a $J^Z$-superconnection and $(A_Z, X_Z)$ is a flat duality superconnection of degree 1.*

**Pf. :** From [9], we have

$$\begin{aligned}
d^M &= d^Z + \nabla^{\Omega(Z)} + i_T, \\
\left(d^M\right)^T &= \left(d^Z\right)^T + \left(\nabla^{\Omega(Z)}\right)^T - T\wedge,
\end{aligned} \tag{222}$$

where in terms of the local framing,

$$\begin{aligned}
d^Z &= E^j \nabla^{TZ}_{e_j}, \\
\nabla^{\Omega(Z)} &= E^\alpha \left(\nabla^{TZ}_{e_\alpha} - \omega_{\alpha j k} E^j I^k\right), \\
i_T &= -\omega_{\alpha\beta j} E^\alpha E^\beta I^j, \\
\left(d^Z\right)^T &= -I^j \nabla^{TZ}_{e_j}, \\
\left(\nabla^{\Omega(Z)}\right)^T &= E^\alpha \left(\nabla^{TZ}_{e_\alpha} - \omega_{\alpha j k} I^j E^k\right), \\
T\wedge &= -\omega_{\alpha\beta j} E^\alpha E^\beta E^j.
\end{aligned} \tag{223}$$

One can check that

$$\begin{aligned}
\left(J^Z\right)^{-1} E^j J^Z &= I^j, \\
\left(J^Z\right)^{-1} I^j J^Z &= E^j, \\
\left(J^Z\right)^{-1} E^\alpha J^Z &= E^\alpha.
\end{aligned} \quad (224)$$

It follows that $A_Z$ is a $J^Z$-superconnection.

As $A_Z - X_Z = d^M$ and $A_Z + X_Z = \left(d^M\right)^T$ are flat superconnections and $A'$ is of total degree 1 with respect to the **Z**-grading coming from $N$, one can reverse the proof of Proposition 26 to show that $(A_Z, X_Z)$ is a flat pair. From (222), we see that $(A_Z, X_Z)$ is of degree 1. It remains to show that $(A_Z, X_Z)$ is $(\cdot,\cdot)_Z$-compatible. From (113), this is equivalent to showing that for $e, e' \in C^\infty(B; \Omega(Z))$,

$$\begin{aligned}
0 &= (d^Z e, e')_Z + (e, d^Z e')_Z, & (225) \\
d(e, e')_Z &= (\nabla^{\Omega(Z)} e, e')_Z + (e, \nabla^{\Omega(Z)} e')_Z, & (226) \\
0 &= (i_T e, e')_Z + (e, i_T e')_Z. & (227)
\end{aligned}$$

Equation (225) follows as in the proof of Lemma 5e. To see (226), recall that if $U$ is a vector field on $B$ then $\nabla_U^{\Omega(Z)}$ is Lie differentiation on $M$ in the direction of the horizontal lift $U^H$ of $U$ [9, Definition 3.2]. As the form $(\cdot,\cdot)_Z$ is diffeomorphism-invariant, (226) follows. Finally, one can check that

$$0 = (I^j e, e')_Z - (e, I^j e')_Z, \quad (228)$$

which implies (227). ∎

Let $\left(E, \nabla^E, (\cdot,\cdot)_E\right)$ be a flat $\epsilon$-symmetric duality bundle on $M$ with an automorphism $J^E$ as before. Put $W = \Omega(M; E)$. Put $(\cdot,\cdot)_W = (\cdot,\cdot)_Z \otimes_{C^\infty(B)} (\cdot,\cdot)_E$, an $\epsilon\epsilon_n$-symmetric nondegenerate bilinear form on $W$. Put $J^W = J^Z \otimes_{C^\infty(B)} J^E$. Let $d^{M,E}$ be exterior differentiation on $W$ and let $\left(d^{M,E}\right)^T$ be its adjoint. Put

$$\begin{aligned}
A &= \frac{1}{2}\left(\left(d^{M,E}\right)^T + d^{M,E}\right), \\
X &= \frac{1}{2}\left(\left(d^{M,E}\right)^T - d^{M,E}\right).
\end{aligned} \quad (229)$$

As in Proposition 30, $A$ is a $J^W$-superconnection and $(A, X)$ is a flat duality superconnection of degree 1. Let $H(Z; E|_Z)$ denote the real vector bundle on $B$ whose fiber over $b \in B$ is isomorphic to $\mathrm{H}(Z_b; E|_{Z_b})$. As in Subsection 3.4, we obtain a flat connection $\nabla^H$ on $H(Z; E|_Z)$, a covariantly-constant $\epsilon\epsilon_n$-symmetric nondegenerate bilinear form $(\cdot, \cdot)_H$ on $H(Z; E|_Z)$ and an automorphism $J^H \in \mathrm{Aut}(H(Z; E|_Z))$.

Define the $J^W$-superconnection $C_t$ and $D_t \in \Omega(B; \mathrm{End}(W))$ as in (201). Explicitly,

$$C_t = \frac{\sqrt{t}}{2}\left(\left(d^{Z,E}\right)^T + d^{Z,E}\right) + \nabla^{W,u} - \frac{1}{2\sqrt{t}} c(T), \tag{230}$$

$$D_t = \frac{\sqrt{t}}{2}\left(\left(d^{Z,E}\right)^T - d^{Z,E}\right) + E^\alpha\left(\omega_{\alpha jk} c^j \widehat{c}^k + \omega_\alpha\left(\nabla^E, J^E\right)\right) - \frac{1}{2\sqrt{t}} \widehat{c}(T).$$

**Remark 5 :** If $J^E$ is covariantly-constant with respect to $\nabla^E$ then $C_{4t}$ is the same as the Bismut superconnection of the twisted vertical signature operator [6, 5].

Define $p(A', J^W(t)) \in \Omega^{4*+1-\epsilon}(B; \mathbf{R})$ and $\widetilde{\eta}(t) \in \Omega^{4*-\epsilon}(B; \mathbf{R})$ by

$$p(A', J^W(t)) = \frac{1}{\sqrt{\epsilon}} \phi \operatorname{tr}\left[J^W e^{-C_t^2}\right],$$

$$\widetilde{\eta}(t) = (2i\pi)^{-1/2} \frac{1}{2t\sqrt{\epsilon}} \phi \operatorname{tr}\left[J^W [N, D_t] e^{-C_t^2}\right]. \tag{231}$$

Let $L\left(TZ, \nabla^{TZ}\right) \in \Omega^{4*}(M)$ denote the Hirzebruch $L$-form of the connection $\nabla^{TZ}$.

**Proposition 31** *The form $p(A', J^W(t))$ is closed. Its de Rham cohomology class $p(A') \in \mathrm{H}^{4*+1-\epsilon}(B; \mathbf{R})$ is independent of $t \in \mathbf{R}^+$. Furthermore,*

$$\frac{\partial}{\partial t} p(A', J^E(t)) = d\widetilde{\eta}(t). \tag{232}$$

*As $t \to \infty$,*

$$p(A', J^W(t)) = p(\nabla^H, J^H) + O(t^{-1/2}),$$
$$\widetilde{\eta}(t) = O(t^{-3/2}). \tag{233}$$

As $t \to 0$,

$$p(A', J^W(t)) = \int_Z L\left(TZ, \nabla^{TZ}\right) \wedge p(\nabla^E, J^E) + O(t),$$
$$\widetilde{\eta}(t) = O(1). \tag{234}$$

**Pf. :** The proofs are similar to those of the analogous statements in [6, 7, 5]. We omit the details. ∎

**Corollary 10** *As elements of* $\mathrm{H}^{4*+1-\epsilon}(B; \mathbf{R})$,

$$p(\nabla^H) = \int_Z L(TZ) \cup p(\nabla^E). \tag{235}$$

**Definition 33** *Define* $\widetilde{\eta}\left(T^H M, g^{TZ}, J^E\right) \in \Omega^{4*-\epsilon}(B)$ *by*

$$\widetilde{\eta}\left(T^H M, g^{TZ}, J^E\right) = -\int_0^\infty \widetilde{\eta}(t) dt. \tag{236}$$

By Proposition 31, $\widetilde{\eta}\left(T^H M, g^{TZ}, J^E\right)$ is well-defined and satisfies

$$d\widetilde{\eta}\left(T^H M, g^{TZ}, J^E\right) = \int_Z L\left(TZ, \nabla^{TZ}\right) \wedge p(\nabla^E, J^E) - p\left(\nabla^H, J^H\right). \tag{237}$$

If $J^E$ is covariantly-constant with respect to $\nabla^E$ then $\widetilde{\eta}\left(T^H M, g^{TZ}, J^E\right)$ is the same as the eta-form of [7] for the twisted vertical signature operator.

## 3.6   Analytic Pushforward of $\overline{L}$-Groups

Let $Z \to M \xrightarrow{\pi} B$ be a fiber bundle as in Subsection 3.5. We assume that $n$ is even. Then if $\epsilon = \pm 1$, $\epsilon \epsilon_n \equiv \epsilon + n \pmod 4$.

**Definition 34** *The pushforward in real cohomology, denoted* $\pi_! : \mathrm{H}^{4*+j}(M; \mathbf{R}) \to \mathrm{H}^{4*+j-n}(B; \mathbf{R})$, *is given by*

$$\pi_!(\tau) = \int_Z L(TZ) \cup \tau. \tag{238}$$

**Definition 35** *The pushforward in $L^0_\epsilon$, denoted $\pi_! : L^0_\epsilon(M) \to L^0_{\epsilon\epsilon_n}(B)$, is generated by*

$$\pi_!\left(\left[\left(E, \nabla^E, (\cdot, \cdot)_E\right)\right]\right) = \left[\left(H(Z; E|_Z), \nabla^H, (\cdot, \cdot)_H\right)\right]. \tag{239}$$

To see that $\pi_!$ is well-defined on $L^0_\epsilon(M)$, suppose first that we are given a short exact sequence (142) on $M$. Using the fact that it splits, we obtain a short exact sequence

$$0 \longrightarrow \mathcal{H}^1 \longrightarrow \mathcal{H}^2 \longrightarrow \mathcal{H}^3 \longrightarrow 0, \tag{240}$$

where $\mathcal{H}^i = \left(H(Z; E^i|_Z), \nabla^{H(Z;E^i|_Z)}, (\cdot, \cdot)_{H(Z;E^i|_Z)}\right)$. Next, if we have an $L^0_\epsilon(M)$-relation of the form $\mathcal{E} = \mathcal{V} \oplus \mathcal{V}^T$ then

$$\begin{aligned}
H(Z; E|_Z) &= H(Z; V|_Z) \oplus H(Z; V|_Z)^T, \\
\nabla^{H(Z;E|_Z)} &= \nabla^{H(Z;V|_Z)} \oplus \left(\nabla^{H(Z;V|_Z)}\right)^T, \\
(\cdot, \cdot)_{H(Z;E|_Z)} &= (\cdot, \cdot)_{H(Z;V|_Z) \oplus H(Z;V|_Z)^T}.
\end{aligned} \tag{241}$$

It follows that $\pi_!$ is well-defined.

It follows from Corollary 10 that there is a commutative diagram

$$\begin{array}{ccc}
L^0_\epsilon(M) & \xrightarrow{p} & \mathrm{H}^{4*+1-\epsilon}(M; \mathbf{R}) \\
\pi_! \downarrow & & \pi_! \downarrow \\
L^0_{\epsilon\epsilon_n}(B) & \xrightarrow{p} & \mathrm{H}^{4*+1-\epsilon-n}(B; \mathbf{R}).
\end{array} \tag{242}$$

Pick a horizontal distribution $T^H M$ and a vertical Riemannian metric $g^{TZ}$ on the fiber bundle.

**Definition 36** *The pushforward in $\widehat{L}^0_\epsilon$, denoted $\pi_! : \widehat{L}^0_\epsilon(M) \to \widehat{L}^0_{\epsilon\epsilon_n}(B)$, is generated by*

$$\pi_!\left(\left[\left(E, J^E, \rho\right)\right]\right) = \left[\left(H(Z; E|_Z), J^H, \int_Z L(TZ, \nabla^{TZ}) \wedge \rho - \widetilde{\eta}\left(T^H M, g^{TZ}, J^E\right)\right)\right]. \tag{243}$$

**Proposition 32** *The pushforward in $\widehat{L}^0_\epsilon$ is well-defined.*

**Pf.** : First, suppose that we have an $\widehat{L}_\epsilon^0(M)$-relation as in (150). For $j \in \{1, 2, 3\}$, put $\widetilde{\eta}^j = \widetilde{\eta}\left(T^H M, g^{TZ}, J^{E^j}\right)$. Define $\sigma \in \Omega^{4*-\epsilon}(B)/\mathrm{im}(d)$ by

$$\begin{aligned}
\sigma &= \widetilde{\eta}^2 - \widetilde{\eta}^1 - \widetilde{\eta}^3 + \widetilde{p}\left(\nabla^{H^2}, J^{H^2}, J^{H^1} \oplus J^{H^3}\right) \\
&\quad - \int_Z L\left(TZ, \nabla^{TZ}\right) \wedge \widetilde{p}\left(\nabla^{E^2}, J^{E^2}, J^{E^1} \oplus J^{E^3}\right).
\end{aligned} \tag{244}$$

Then we must show that $\sigma = 0$. If $J^{E^2} = J^{E^1} \oplus J^{E^3}$ then we are in a direct sum situation and $\sigma = 0$. Thus it suffices to show that if $J_s^{E^2}$ is a smooth 1-parameter family of $J^{E^2}$'s then $\partial \sigma / \partial s = 0$. Equivalently, it suffices to show that if $\left(E, \nabla^E, (\cdot, \cdot)_E\right)$ is an $L_\epsilon^0(M)$-generator and $J_s^E$ is a smooth 1-parameter family of $J^E$'s parameterized by $s \in [0, 1]$ then

$$\begin{aligned}
\frac{\partial \widetilde{\eta}\left(T^H M, g^{TZ}, J_s^E\right)}{\partial s} &= \int_Z L\left(TZ, \nabla^{TZ}\right) \wedge \\
&\quad (2i\pi)^{-1/2} \frac{1}{4\sqrt{\epsilon}} \phi \, \mathrm{tr}\left[J^E \left[(J^E)^{-1} \frac{\partial J^E}{\partial s}, \omega_E\right] e^{\frac{\omega_E^2}{4}}\right] \\
&\quad - (2i\pi)^{-1/2} \frac{1}{4\sqrt{\epsilon}} \phi \, \mathrm{tr}\left[J^H \left[(J^H)^{-1} \frac{\partial J^H}{\partial s}, \omega_H\right] e^{\frac{\omega_H^2}{4}}\right].
\end{aligned} \tag{245}$$

Define $\widetilde{M}, \widetilde{B}, T^H \widetilde{M}$ and $g^{T\widetilde{Z}}$ as in the proof of Lemma 1. Put $\widetilde{E} = \pi_M^* E$, $\nabla^{\widetilde{E}} = \pi_M^* \nabla^E$ and $(\cdot, \cdot)_{\widetilde{E}} = \pi_M^* (\cdot, \cdot)_E$. We abbreviate

$$\left(H\left(\widetilde{Z}; \widetilde{E}|_{\widetilde{Z}}\right), \nabla^{H\left(\widetilde{Z}; \widetilde{E}|_{\widetilde{Z}}\right)}, (\cdot, \cdot)_{H\left(\widetilde{Z}; \widetilde{E}|_{\widetilde{Z}}\right)}\right)$$

by $\left(\widetilde{H}, \nabla^{\widetilde{H}}, (\cdot, \cdot)_{\widetilde{H}}\right)$. Then

$$\left(\widetilde{H}, \nabla^{\widetilde{H}}, (\cdot, \cdot)_{\widetilde{H}}\right) = \pi_B^* \left(H, \nabla^H, (\cdot, \cdot)_H\right). \tag{246}$$

Define $J^{\widetilde{E}} \in \mathrm{Aut}(\widetilde{E})$ so that its restriction to $\{s\} \times M$ is $J_s^E$.

The exterior differentiation on $\widetilde{B}$ is given by (36). Consider the eta-form $\widetilde{\eta}\left(T^H \widetilde{M}, g^{T\widetilde{Z}}, J_s^{\widetilde{E}}\right) \in \Omega^{4*-\epsilon}(\widetilde{B})$. By (237),

$$d\widetilde{\eta}\left(T^H \widetilde{M}, g^{T\widetilde{Z}}, J_s^{\widetilde{E}}\right) = \int_{\widetilde{Z}} L\left(T\widetilde{Z}, \nabla^{T\widetilde{Z}}\right) \wedge p(\nabla^{\widetilde{E}}, J^{\widetilde{E}}) - p\left(\nabla^{\widetilde{H}}, J^{\widetilde{H}}\right). \tag{247}$$

By construction,
$$L\left(T\tilde{Z}, \nabla^{T\tilde{Z}}\right) = \pi_M^* L\left(TZ, \nabla^{TZ}\right). \tag{248}$$

Equations (36) and (247) give that modulo im($d$),
$$\frac{\partial \tilde{\eta}}{\partial s} = \int_Z L\left(TZ, \nabla^{TZ}\right) \wedge i_{\partial_s} p\left(\nabla^{\tilde{E}}, J^{\tilde{E}}\right) - i_{\partial_s} p\left(\nabla^{\tilde{H}}, J^{\tilde{H}}\right). \tag{249}$$

Equation (245) now follows.

Next, supppose that we have a $\widehat{L}_\epsilon^0(M)$-relation $\mathcal{V} \oplus \mathcal{V}^T = (V \oplus V^T, J_0^{V \oplus V^T}, 0)$. Then with respect to the splitting
$$W = \Omega(M; V) \oplus \Omega(M; V^T), \tag{250}$$

one can see that $J^W[N, D_t] e^{-C_t^2}$ is off-diagonal. Thus $\tilde{\eta}\left(T^H M, g^{TZ}, J_0^{V \oplus V^T}\right)$ vanishes. Furthermore, it follows from Definition 30 that
$$J^{H(Z;E|_Z)} = J_0^{H(Z;V|_Z) \oplus H(Z;V|_Z)^T}. \tag{251}$$

Thus $\pi_!\left(\mathcal{V} \oplus \mathcal{V}^T\right)$ is an $\widehat{L}_{\epsilon\epsilon_n}^0(B)$-relation. The proposition follows. ∎

**Proposition 33** *The pushforward in $\widehat{L}_\epsilon^0$ restricts to a pushforward*
$$\pi_! : \overline{L}_\epsilon^0(M) \to \overline{L}_{\epsilon\epsilon_n}^0(B). \tag{252}$$

**Pf.** : It is enough to show that there is a commutative diagram
$$\begin{array}{ccc} \widehat{L}_\epsilon^0(M) & \xrightarrow{p'} & \Omega^{4*+1-\epsilon}(M) \\ \pi_! \downarrow & & \pi_! \downarrow \\ \widehat{L}_{\epsilon\epsilon_n}^0(B) & \xrightarrow{p'} & \Omega^{4*+1-\epsilon-n}(B). \end{array} \tag{253}$$

This follows from (237). ∎

**Proposition 34** *The pushforward in $\overline{L}_\epsilon^0$ is independent of the choices of $T^H M$ and $g^{TZ}$.*

**Pf.** : We use the notation of the proof of Proposition 9. Given a finite set of generators $\mathcal{E}^j = \left(E^j, J^{E^j}, \rho^j\right)$ in $\widehat{L}^0_\epsilon(M)$, put $\widetilde{\mathcal{E}^j} = \pi_M^* \mathcal{E}^j$. If $\sum_j n_j [\mathcal{E}^j]$ lies in $\overline{L}^0_\epsilon(M)$ then $\sum_j n_j \left[\widetilde{\mathcal{E}^j}\right]$ lies in $\overline{L}^0_\epsilon(\widetilde{M})$ and

$$\widetilde{l} = \widetilde{\pi}_! \sum_j n_j \left[\widetilde{\mathcal{E}^j}\right] \tag{254}$$

lies in $\overline{L}^0_{\epsilon\epsilon_n}(\widetilde{B})$. By construction, $i_0^* \widetilde{l}$ is the pushforward of $\sum_j n_j [\mathcal{E}^j]$ using $(T^H M, g^{TZ})$ and $i_1^* \widetilde{l}$ is the pushforward of $\sum_j n_j [\mathcal{E}^j]$ using $(T'^H M, g'^{TZ})$. The proposition now follows from Proposition 22. ∎

# A  Results from [9]

In this appendix we describe results from [9] on flat complex vector bundles and their direct images.

## A.1  Characteristic Classes of Flat Complex Vector Bundles

Let $B$ be a smooth connected compact manifold. If $E$ is a complex vector bundle over $B$, we let $C^\infty(B; E)$ denote the smooth sections of $E$. We let $\Lambda(T^*B)$ denote the complexified exterior bundle of $B$ and $\Omega(B)$ denote the smooth sections of $\Lambda(T^*B)$. We put $\Omega(B; E) = C^\infty(B; \Lambda(T^*B) \otimes E)$. We say that a differential form is real if it can be written with real coefficients.

Let $\phi : \Omega(B) \to \Omega(B)$ be the linear map such that for all homogeneous $\omega \in \Omega(B)$,

$$\phi\, \omega = (2i\pi)^{-(\deg \omega)/2} \omega. \tag{255}$$

Let $E$ be a complex vector bundle on $B$, endowed with a flat connection $\nabla^E$. Let $h^E$ be a positive-definite Hermitian metric on $E$. We do not require that $\nabla^E$ be compatible with $h^E$. Define $\omega(\nabla^E, h^E) \in \Omega^1(B; \mathrm{End}(E))$ by

$$\omega(\nabla^E, h^E) = (h^E)^{-1}(\nabla^E h^E). \tag{256}$$

With respect to a locally-defined covariantly-constant basis of $E$, $h^E$ is locally a Hermitian matrix-valued function on $B$ and we can write $\omega(\nabla^E, h^E)$ more

simply as
$$\omega(\nabla^E, h^E) = (h^E)^{-1} dh^E. \tag{257}$$

**Definition 37** *For $k$ a positive odd integer, define $c_k(\nabla^E, h^E) \in \Omega^k(B)$ by*
$$c_k(\nabla^E, h^E) = (2i\pi)^{-\frac{k-1}{2}} \, 2^{-k} \operatorname{Tr}\left[\omega^k(\nabla^E, h^E)\right]. \tag{258}$$

*Let $c(\nabla^E, h^E) \in \Omega^{odd}(B)$ be the formal sum*
$$c(\nabla^E, h^E) = \sum_{j=0}^{\infty} \frac{1}{j!} \, c_{2j+1}(\nabla^E, h^E). \tag{259}$$

**Proposition 35** *[9] The form $c_k(\nabla^E, h^E)$ is real and closed. Its de Rham cohomology class is independent of $h^E$.*

**Definition 38** *We will denote the de Rham cohomology class of $c_k(\nabla^E, h^E)$ by $c_k(\nabla^E) \in \mathrm{H}^k(B; \mathbf{R})$ and the de Rham cohomology class of $c(\nabla^E, h^E)$ by $c(\nabla^E) \in \mathrm{H}^{odd}(B; \mathbf{R})$.*

The classes $c_k(\nabla^E)$ are the characteristic classes (of flat vector bundles) which are of interest to us. A more topological description of them can be given as follows. Let $V$ be a finite-dimensional complex vector space. Let $\mathrm{H}^*_c(GL(V); \mathbf{R})$ denote the continuous group cohomology of $GL(V)$, meaning the cohomology of the complex of Eilenberg-Maclane cochains on $GL(V)$ which are continous in their arguments. Let $GL(V)_\delta$ denote $GL(V)$ with the discrete topology and let $BGL(V)_\delta$ denote its classifying space. The cohomology group $\mathrm{H}^*(BGL(V)_\delta; \mathbf{R})$ is isomorphic to the (discrete) group cohomology $\mathrm{H}^*(GL(V); \mathbf{R})$. There is a forgetful map
$$\mu_V : \mathrm{H}^*_c(GL(V); \mathbf{R}) \to \mathrm{H}^*(BGL(V)_\delta; \mathbf{R}). \tag{260}$$

Fix a basepoint $* \in B$. Put $\Gamma = \pi_1(B, *)$ and let $h : B \to B\Gamma$ be the classifying map for the universal cover of $B$, defined up to homotopy. Let $V$ be the fiber of $E$ above $*$. The holonomy of $E$ is a homomorphism $r : \Gamma \to GL(V)$, and induces a map $Br : B\Gamma \to BGL(V)_\delta$. Then the flat bundle $E$ is classified by the homotopy class of maps $\nu = Br \circ h : B \to BGL(V)_\delta$. One can show that there is a class $c_{k,V} \in \mathrm{H}^k(GL(V); \mathbf{R})$ such that $c_k(\nabla^E) = \nu^*(c_{k,V})$, and a class $C_{k,V} \in \mathrm{H}^k_c(GL(V); \mathbf{R})$ such that $c_{k,V} = \mu_V(C_{k,V})$. For example, $C_{1,V}$ is given by the homomorphism $g \to \ln|\det(g)|$ from $GL(V)$ to $(\mathbf{R}, +)$.

Put $G = GL(V)$ and $K = U(V)$. Denote the Lie algebras of $G$ and $K$ by $\gamma = gl(V)$ and $\kappa = u(V)$, respectively. The quotient space $\gamma/\kappa$ is isomorphic to the space of Hermitian endomorphisms of $V$, and carries an adjoint representation of $K$. One has that $\mathrm{H}_c^*(GL(V); \mathbf{R})$ is isomorphic to $\mathrm{H}^*(\gamma, K; \mathbf{R})$, the cohomology of the complex $\mathrm{C}^*(\gamma, K; \mathbf{R}) = \mathrm{Hom}_K\left(\Lambda^*(\gamma/\kappa), \mathbf{R}\right)$ [11, Chapter IX, §5]. In fact, the differential of this complex vanishes, and so $\mathrm{H}_c^*(GL(V); \mathbf{R}) = \mathrm{C}^*(\gamma, K; \mathbf{R})$ [11, Chapter II, Corollary 3.2]. Thus the classes $\{c_k(\nabla^E)\}$ arise indirectly from $K$-invariant forms on $\gamma/\kappa$. It is possible to see the relationship between $c_k(\nabla^E)$ and $\mathrm{C}^k(\gamma, K; \mathbf{R})$ more directly [9, §1g]. In particular, for $k$ odd define a $k$-form $\Phi_k$ on $\gamma/\kappa$ by sending Hermitian endomorphisms $M_1, \ldots, M_k$ to

$$\Phi_k(M_1, \ldots, M_k) = \sum_{\sigma \in S_k} (-1)^{\mathrm{sign}(\sigma)} \, \mathrm{Tr}\left[M_{\sigma(1)} \ldots M_{\sigma(k)}\right]. \tag{261}$$

Then $\Phi_k$ is an element of $\mathrm{C}^k(\gamma, K; \mathbf{R})$ which, up to an overall multiplicative constant, corresponds to $C_{k,V}$.

The compact dual of the symmetric space $G/K$ is $G^d/K$, where $G^d = U(V) \times U(V)$. Let $\gamma^d = u(V) \oplus u(V)$ be the Lie algebra of $G^d$. Duality gives an isomorphism between $\mathrm{H}^*(\gamma, K; \mathbf{R})$ and $\mathrm{H}^*(\gamma^d, K; \mathbf{R}) = \mathrm{H}^*(U(V); \mathbf{R}) = \Lambda(x_1, x_3, \ldots, x_{2\dim(V)-1})$. It follows that the classes $\{C_{2j-1,V}\}_{j=1}^{\dim(V)}$ are algebraically independent.

## A.2 The Superconnection Formalism

For background information on superconnections we refer to [6, 5, 30]. Let $E = E_+ \oplus E_-$ be a $\mathbf{Z}_2$-graded finite-dimensional complex vector bundle on $B$. Let $\tau$ be the involution of $E$ defining the $\mathbf{Z}_2$-grading, so that $\tau|_{E_\pm} = \pm I$. Then $\mathrm{End}(E)$ is a $\mathbf{Z}_2$-graded bundle of algebras over $B$, whose even (resp. odd) elements commute (resp. anticommute) with $\tau$. Given $a \in C^\infty(B; \mathrm{End}(E))$, we define its supertrace $\mathrm{Tr}_s[a] \in C^\infty(B)$ by

$$\mathrm{Tr}_s[a] = \mathrm{Tr}[\tau a]. \tag{262}$$

Given $\omega \in \Omega(B)$ and $a \in C^\infty(B; \mathrm{End}(E))$, put

$$\mathrm{Tr}_s[\omega \cdot a] = \omega \mathrm{Tr}_s[a]. \tag{263}$$

Then $\mathrm{Tr}_s$ extends to a linear map from $\Omega(B; \mathrm{End}(E))$ to $\Omega(B)$.

Given $\alpha, \alpha' \in \Omega(B; \mathrm{End}(E))$, we define their supercommutator $[\alpha, \alpha'] \in \Omega(B; \mathrm{End}(E))$ by

$$[\alpha, \alpha'] = \alpha \alpha' - (-1)^{(\deg \alpha)(\deg \alpha')} \alpha' \alpha. \tag{264}$$

A basic fact is that $\mathrm{Tr}_s$ vanishes on supercommutators [30].

Let $\nabla^E$ be a connection on $E$ which preserves the splitting $E = E_+ \oplus E_-$. Then $\nabla^E$ decomposes as $\nabla^E = \nabla^{E_+} \oplus \nabla^{E_-}$. Let $S$ be an odd element of $\Omega(B; \mathrm{End}(E))$. By definition, $\nabla^E + S$ gives a superconnection $A$ on $E$. That is, there is a $\mathbf{C}$-linear map

$$A : C^\infty(B; E) \to \Omega(B; E) \tag{265}$$

which is odd with respect to the total $\mathbf{Z}_2$-gradings and satisfies the Leibniz rule. We can extend $A$ to an odd $\mathbf{C}$-linear endomorphism of $\Omega(B; E)$. By definition, the curvature of $A$ is $A^2$, an even $C^\infty(B)$-linear endomorphism of $\Omega(B; E)$ which is given by multiplication by an even element of $\Omega(B; \mathrm{End}(E))$.

In what follows, we will say that a holomorphic function $f : \mathbf{C} \to \mathbf{C}$ is real if for all $a \in \mathbf{C}$, we have $f(\overline{a}) = \overline{f(a)}$.

## A.3 Characteristic Classes and Torsion Forms of Flat Superconnections

**Definition 39** *A superconnection $A'$ on $E$ is flat if its curvature vanishes, i.e. if $A'^2 = 0$.*

Hereafter we assume that $A'$ is flat. Let $h^E$ be a Hermitian metric on $E$ such that $E_+$ and $E_-$ are orthogonal. Then there is a flat superconnection $A'^*$ on $E$ which is the adjoint of $A'$ with respect to $h^E$. Define an odd element of $\Omega(B; \mathrm{End}(E))$ by

$$X = \frac{1}{2}(A'^* - A'). \tag{266}$$

**Definition 40** *Let $f : \mathbf{C} \to \mathbf{C}$ be a holomorphic real odd function. Put*

$$f\left(A', h^E\right) = (2i\pi)^{1/2} \phi \, \mathrm{Tr}_s[f(X)] \in \Omega(B). \tag{267}$$

**Proposition 36** *[9] The differential form $f\left(A', h^E\right)$ is real, odd and closed. Its de Rham cohomology class is independent of $h^E$.*

We will denote the de Rham cohomology class of $f\left(A', h^E\right)$ by $f(A') \in \mathrm{H}^{odd}(B; \mathbf{R})$.

Suppose now that $E = \bigoplus_{i=0}^n E^i$ is a **Z**-graded complex vector bundle on $B$. Put

$$E_+ = \bigoplus_{i \text{ even}} E^i \ , \quad E_- = \bigoplus_{i \text{ odd}} E^i. \tag{268}$$

Then $E = E_+ \oplus E_-$ is a $\mathbf{Z}_2$-graded vector bundle, to which we may apply the above formalism.

Let $A'$ be a superconnection on $E$. We can expand $A$ as

$$A = \sum_{j \geq 0} A_j, \tag{269}$$

where $A_j$ is of partial degree $j$ with respect to the **Z**-grading on $\Lambda(T^*B)$.

**Definition 41** *We say that $A'$ is of total degree 1 if*
- *$A'_1$ is a connection on $E$ which preserves the **Z**-grading.*
- *For $j \in \mathbf{N} - \{1\}$, $A'_j$ is an element of $\Omega^j(B; \mathrm{Hom}(E^\bullet, E^{\bullet+1-j}))$.*

In what follows, we will assume that $A'$ is a flat superconnection of total degree 1. Put

$$v = A'_0 \ , \quad \nabla^E = A'_1. \tag{270}$$

Clearly $v \in C^\infty(B; \mathrm{Hom}(E^\bullet, E^{\bullet+1}))$. The flatness of $A'$ implies that

$$v^2 = \left[\nabla^E, v\right] = \left(\nabla^E\right)^2 + [v, A'_2] = 0. \tag{271}$$

As $v^2 = 0$, we have a cochain complex of vector bundles

$$(E, v): 0 \longrightarrow E^0 \xrightarrow{v} E^1 \xrightarrow{v} \cdots \xrightarrow{v} E^n \longrightarrow 0. \tag{272}$$

**Definition 42** *For $b \in B$, let $H(E, v)_b = \bigoplus_{i=0}^n \mathrm{H}^i(E, v)_b$ be the cohomology of the complex $(E, v)_b$.*

Using (271), one can show that there a **Z**-graded complex vector bundle $H(E,v)$ on $B$ whose fiber over $b \in B$ is $H(E,v)_b$, and a natural flat connection $\nabla^{H(E,v)}$ on $H(E,v)$.

Let $h^E$ be a Hermitian metric on $E$ such that the $E^i$'s are mutually orthogonal. Put $A'' = A'^*$, the adjoint superconnection to $A'$ with respect to $h^E$. Let $v^* \in C^\infty(B; \text{Hom}(E^\bullet, E^{\bullet-1}))$ be the adjoint of $v$ with respect to $h^E$. From finite-dimensional Hodge theory, there is an isomorphism

$$H(E,v) \cong \text{Ker}(v^*v + vv^*). \tag{273}$$

Being a subbundle of $E$, the vector bundle $\text{Ker}(v^*v+vv^*)$ inherits a Hermitian metric $h^{Ker}$ from the Hermitian metric $h^E$ on $E$. Let $h^{H(E,v)}$ denote the Hermitian metric on $H(E,v)$ obtained from $h^{Ker}$ via the isomorphism (273).

Let $N \in \text{End}(E)$ be the number operator of $E$, i.e. $N$ acts on $E^i$ by multiplication by $i$. Extend $N$ to an element of $C^\infty(B; \text{End}(E))$.

**Definition 43** *For $t > 0$, let $C'_t$ be the flat superconnection on $E$ of total degree 1 given by*

$$C'_t = t^{N/2} \; A' \; t^{-N/2}. \tag{274}$$

*and let $C''_t$ be the flat superconnection on $E$ given by*

$$C''_t = t^{-N/2} \; A'' \; t^{N/2}. \tag{275}$$

The superconnections $C'_t$ and $C''_t$ are adjoint with respect to $h^E$. We have

$$\begin{aligned} C'_t &= \sum_{j \geq 0} t^{(1-j)/2} A'_j, \\ C''_t &= \sum_{j \geq 0} t^{(1-j)/2} A''_j. \end{aligned} \tag{276}$$

Define an odd element of $\Omega(B; \text{End}(E))$ by

$$D_t = \frac{1}{2}\left(C''_t - C'_t\right). \tag{277}$$

**Definition 44** *Define a real even differential form on $B$ by*

$$f^\wedge\left(C'_t, h^E\right) = \phi \, \text{Tr}_s\left[\frac{N}{2} f'(D_t)\right] \in \Omega(B). \tag{278}$$

**Proposition 37** *[9]* One has

$$\frac{\partial}{\partial t} f\left(C'_t, h^E\right) = \frac{1}{t} df^\wedge \left(C'_t, h^E\right). \tag{279}$$

Let $d(H(E,v))$ be the constant integer-valued function on $B$

$$d(H(E,v)) = \sum_{i=0}^{n} (-1)^i \, i \, \mathrm{rk}\left(H^i(E,v)\right). \tag{280}$$

Hereafter we take $f(z) = z \, \exp(z^2)$.

**Proposition 38** *[9]* As $t \to +\infty$,

$$\begin{aligned} f\left(C'_t, h^E\right) &= f\left(\nabla^{H(E,v)}, h^{H(E,v)}\right) + \mathcal{O}\left(\frac{1}{\sqrt{t}}\right) \\ f^\wedge \left(C'_t, h^E\right) &= d(H(E,v)) \frac{f'(0)}{2} + \mathcal{O}\left(\frac{1}{\sqrt{t}}\right). \end{aligned} \tag{281}$$

Now consider the special case when the vector bundle $E$ has not only a flat superconnection, but has a flat connection. Let

$$(E,v) : 0 \longrightarrow E^0 \xrightarrow{v} E^1 \xrightarrow{v} \cdots \xrightarrow{v} E^n \longrightarrow 0 \tag{282}$$

be a flat complex of complex vector bundles. That is,

$$\nabla^E = \bigoplus_{i=0}^{n} \nabla^{E^i} \tag{283}$$

is a flat connection on $E = \bigoplus_{i=0}^{n} E^i$ and $v$ is a flat cochain map, meaning

$$\left(\nabla^E\right)^2 = 0 \;,\; v^2 = 0 \;,\; \left[\nabla^E, v\right] = 0. \tag{284}$$

Take $A' = v + \nabla^E$. For $t > 0$, put

$$\begin{aligned} C'_t &= \sqrt{t}\, v + \nabla^E, \\ C''_t &= \sqrt{t}\, v^* + \left(\nabla^E\right)^*. \end{aligned} \tag{285}$$

Then $C'_t$ is a flat superconnection of total degree 1. Let $d(E)$ be the constant integer-valued function on $B$ given by

$$d(E) = \sum_{i=0}^{n} (-1)^i \, i \, \mathrm{rk}\left(E^i\right). \tag{286}$$

**Proposition 39** *[9] As $t \to 0$,*

$$\begin{aligned} f\left(C'_t, h^E\right) &= f\left(\nabla^E, h^E\right) + \mathcal{O}(t), \\ f^\wedge\left(C'_t, h^E\right) &= d(E)\frac{f'(0)}{2} + \mathcal{O}(t). \end{aligned} \quad (287)$$

**Proposition 40** *[9] As elements of $\mathrm{H}^{odd}(B; \mathbf{R})$,*

$$f\left(\nabla^E\right) = f\left(\nabla^{H(E,v)}\right). \quad (288)$$

We now refine (40) to a statement about differential forms on $B$.

**Definition 45** *Define a real even differential form on $B$ by*

$$\begin{aligned} T_f(A', h^E) &= -\int_0^{+\infty} \left[ f^\wedge\left(C'_t, h^E\right) - d(H(E,v))\frac{f'(0)}{2} \right. \\ &\quad \left. - [d(E) - d(H(E,v))]\frac{f'(\frac{i\sqrt{t}}{2})}{2} \right] \frac{dt}{t}. \end{aligned} \quad (289)$$

**Remark 6 :** By Propositions 38 and 39, the integrand in (289) is integrable. We will call $T_f(A', h^E)$ a torsion form.

**Proposition 41** *[9] One has*

$$dT_f(A', h^E) = f\left(\nabla^E, h^E\right) - f\left(\nabla^{H(E,v)}, h^{H(E,v)}\right). \quad (290)$$

**Proof :** This follows from Propositions 37, 38 and 39. ∎

Upon passing to de Rham cohomology, Proposition 41 implies Proposition 40. Up to an overall multiplicative constant, $T_0(A', h^E)$ is the function which to a point $b \in B$ assigns the torsion of the cochain complex $(E, v)_b$, in the sense of [27, 31].

We now relate the constructions of this subsection to those of Subsection A.1. As in Subsection A.1, let $E$ be a complex vector bundle on $B$, endowed with a flat connection $\nabla^E$. We can consider $E$ to be a $\mathbf{Z}_2$-graded vector bundle with $E_+ = E_\mathbf{C}$ and $E_- = 0$. Let $h^E$ be a Hermitian metric on $E$. Taking the flat superconnection $A'$ to be $\nabla^E$, the $X$ of (266) is given by

$$X = \frac{\omega\left(\nabla^E, h^E\right)}{2}. \tag{291}$$

From Definition 40, we have

$$f\left(\nabla^E, h^E\right) = c\left(\nabla^E, h^E\right). \tag{292}$$

Suppose now that we have an exact sequence of flat complex vector bundles on $B$ :

$$0 \longrightarrow E^1 \xrightarrow{v} E^2 \xrightarrow{v} \ldots \xrightarrow{v} E^n \longrightarrow 0. \tag{293}$$

Let $\nabla^{E^i}$ be the flat connection on $E^i$ and let $h^{E^i}$ be a positive-definite Hermitian metric on $E^i$. Then the **Z**-graded vector bundle $E = \oplus_{i=1}^n E^i$ acquires a flat connection $\nabla^E = \oplus_{i=1}^n \nabla^{E^i}$ and a Hermitian metric $h^E = \oplus_{i=1}^n h^{E^i}$. Let $A'$ be the flat superconnection on $E$ given by $A' = v + \nabla^E$.

**Definition 46** *In this special case, define $T_f(A', h^E) \in \Omega^{even}(B)/\mathrm{im}(d)$ as in (289), taken modulo $\mathrm{im}(d)$.*

Then from equation (41),

$$dT_f(A', h^E) = \sum_{i=1}^n (-1)^i c(\nabla^{E^i}, h^{E^i}). \tag{294}$$

We can think of $T_f(A', h^E)$ as an analog of the Bott-Chern class [12].

## A.4 Fiber Bundles

Let $Z \to M \xrightarrow{\pi} B$ be a smooth fiber bundle with connected compact base $B$ and connected closed fibers $Z_b = \pi^{-1}(b)$. Let $F$ be a flat complex vector bundle on $M$. Let $H(Z; F|_Z)$ denote the **Z**-graded complex vector bundle on $B$ whose fiber over $b \in B$ is isomorphic to the cohomology group $\mathrm{H}^*(Z_b, F|_{Z_b})$. It has a canonical flat connection $\nabla^{H(Z;F|_Z)}$ which preserves the **Z**-grading. Let $TZ$ be the vertical tangent bundle of the fiber bundle and let $o(TZ)$ be its orientation bundle, a flat real line bundle on $M$. Let $e(TZ) \in \mathrm{H}^{\dim(Z)}(M; o(TZ))$ be the Euler class of $TZ$. Put $f(z) = z \exp(z^2)$.

**Proposition 42** *[9] One has an equality in $\mathrm{H}^{odd}(B; \mathbf{R})$:*

$$f\left(\nabla^{H(Z;F|_Z)}\right) = \int_Z e(TZ) \cup f(\nabla^F). \tag{295}$$

In fact, one can refine (295) to a statement about differential forms on $B$. First, equip the fiber bundle with a horizontal distribution $T^H M$. Let $W$ be the infinite-dimensional **Z**-graded vector bundle on $B$ whose fiber over $b \in B$ is isomorphic to $\Omega\left(Z_b; F\big|_{Z_b}\right)$. Then

$$C^\infty(B; W) \simeq C^\infty(M; \Lambda(T^*Z) \otimes F) \tag{296}$$

and there is an isomorphism of **Z**-graded vector spaces

$$\Omega(M; F) \simeq \Omega(B; W). \tag{297}$$

Let $N$ be the number operator of $W$; it acts as multiplication by $i$ on $C^\infty(M; \Lambda^i(T^*Z) \otimes F)$.

The exterior differentiation operator $d^M$, acting on $\Omega(M; F)$, defines a flat superconnection on $W$ of total degree 1. In terms of the **Z**-grading on $\Lambda(T^*B)$, $d^M$ can be decomposed as

$$d^M = d^Z + \nabla^W + i_T, \tag{298}$$

where $d^Z$ is vertical exterior differentiation, $\nabla^W$ is a natural connection on $W$ which preserves the **Z**-grading and $i_T$ is interior multiplication by the curvature $T$ of the fiber bundle, a $TZ$-valued horizontal 2-form on $M$. For $t > 0$, put

$$\begin{aligned} C'_t &= t^{N/2} \, d^M \, t^{-N/2} \\ &= \sqrt{t}\, d^Z + \nabla^W + \frac{1}{\sqrt{t}}\, i_T. \end{aligned} \tag{299}$$

Now equip the fiber bundle with a vertical Riemannian metric $g^{TZ}$ and the flat vector bundle $F$ with a Hermitian metric $h^F$. Then $W$ acquires an $L^2$-inner product $h^W$. There is a canonical metric-compatible connection $\nabla^{TZ}$ on $TZ$ [6, 5]. The vector bundle $H(Z; F|_Z)$ acquires a Hermitian metric $h^{H(Z;F|_Z)}$ from Hodge theory. Let $C''_t$ be the adjoint superconnection to $C'_t$ with respect to $h^W$. That is,

$$C''_t = \sqrt{t}\left(d^Z\right)^* + \left(\nabla^W\right)^* - \frac{1}{\sqrt{t}}\, (T \wedge). \tag{300}$$

Define $D_t$, an odd element of $\Omega(B; \text{End}(W))$, by

$$D_t = \frac{1}{2}\left(C_t'' - C_t'\right). \qquad (301)$$

For $t > 0$, define a real odd differential form on $B$ by

$$f\left(C_t', h^W\right) = (2i\pi)^{1/2} \; \phi \, \text{Tr}_s\left[f\left(D_t\right)\right] \qquad (302)$$

and a real even differential form on $B$ by

$$f^\wedge\left(C_t', h^W\right) = \phi \, \text{Tr}_s\left[\frac{N}{2} f'(D_t)\right]. \qquad (303)$$

**Proposition 43** *[9] For any $t > 0$,*

$$\frac{\partial}{\partial t} f\left(C_t', h^W\right) = \frac{1}{t} \, df^\wedge\left(C_t', h^W\right). \qquad (304)$$

Put

$$\chi'(Z; F) = \sum_{i=0}^{\dim(Z)} (-1)^i \, i \, \text{rk}\left(H^i(Z; F|_Z)\right), \qquad (305)$$

an integer-valued constant function on $B$.

**Proposition 44** *[9] As $t \to 0$,*

$$
\begin{aligned}
f\left(C_t', h^W\right) &= \int_Z e\left(TZ, \nabla^{TZ}\right) f\left(\nabla^F, h^F\right) + \mathcal{O}(t) \text{ if } \dim(Z) \text{ is even} \\
&= \mathcal{O}\left(\sqrt{t}\right) \text{ if } \dim(Z) \text{ is odd}, \\
f^\wedge\left(C_t', h^W\right) &= \frac{1}{4}\dim(Z)\text{rk}(F)\chi(Z) + \mathcal{O}(t) \text{ if } \dim(Z) \text{ is even}, \\
&= \mathcal{O}\left(\sqrt{t}\right) \text{ if } \dim(Z) \text{ is odd}. \qquad (306)
\end{aligned}
$$

As $t \to +\infty$,

$$
\begin{aligned}
f\left(C_t', h^W\right) &= f\left(\nabla^{H(Z;F|_Z)}, h^{H(Z;F|_Z)}\right) + \mathcal{O}\left(\frac{1}{\sqrt{t}}\right), \\
f^\wedge\left(C_t', h^W\right) &= \frac{\chi'(Z; F)}{2} + \mathcal{O}\left(\frac{1}{\sqrt{t}}\right). \qquad (307)
\end{aligned}
$$

**Definition 47** *The analytic torsion form $\mathcal{T}\left(T^H M, g^{TZ}, h^F\right)$, a real even differential form on B, is given by*

$$\mathcal{T}\left(T^H M, g^{TZ}, h^F\right) = -\int_0^{+\infty} \left[f^\wedge\left(C'_t, h^W\right) - \frac{\chi'(Z;F)}{2}f'(0) \right.$$
$$\left. -\left(\frac{\dim(Z)\mathrm{rk}(F)\chi(Z)}{4} - \frac{\chi'(Z;F)}{2}\right)f'(\frac{i\sqrt{t}}{2})\right]\frac{dt}{t}.$$
(308)

**Remark 7 :** It follows from Proposition 44 that the integrand of (308) is integrable.

**Proposition 45** *[9] One has*

$$d\mathcal{T}\left(T^H M, g^{TZ}, h^F\right) = \int_Z e\left(TZ, \nabla^{TZ}\right) \wedge f\left(\nabla^F, h^F\right) - f\left(\nabla^{H(Z;F|_Z)}, h^{H(Z;F|_Z)}\right).$$
(309)

**Proof :** This follows from Propositions 43 and 44. ∎

Upon passing to de Rham cohomology, Proposition 45 implies Proposition 42.

**Remark 8 :** One can extend the results of this subsection to the case when the fiber has boundary, by using the doubling trick of [23, Section IX].